\documentclass[12pt]{article}
\usepackage{amssymb}
\usepackage[cp1251]{inputenc}
\usepackage[english]{babel}
\usepackage{amsmath}

\def\mapr#1{\smash{\mathop{\buildrel{#1}\over\longrightarrow}}}

\newlength{\myunit}
\newcounter{mycount}
\def\R{{\bf R}}
\def\F{{\bf F}}

\def\Pr{{\bf Pr}}

\def\rank{\mbox{rank}}

\def\const{\mbox{\bf const }}

\def\inxx#1{}

\def\MYdef{\mathrel{\stackrel{\rm def}=}}

\newtheorem{thm}{Theorem}
\newtheorem{lm}{Lemma}
\newtheorem{df}{Definition}
\newtheorem{pr}{Proposition}
\newtheorem{cor}{Corollary}
\newtheorem{rem}{Remark}

\newcommand\proof{{\bf Proof. }\nobreak\noindent}

\def\mathrm#1{\hbox{\rm #1 }}

\author{Anwar A. Irmatov}
\title{ Singularity of $\{\pm 1\}$-matrices and asymptotics of the number of threshold functions}

\date{}

\begin{document}

\renewcommand\refname{\centering  \sc References }

\maketitle
\begin{abstract}
Two results concerning the number of threshold functions $P(2, n)$ and the probability ${\mathbb P}_n$ that a random $n\times n$ Bernoulli matrix is singular are established. We introduce a supermodular function $\eta^{\bigstar}_n : 2^{{\bf RP}^n}_{fin} \to \mathbb{Z}_{\geq 0},$ defined on finite subsets of ${\bf RP}^n,$ that allows us to obtain a lower bound for $P(2, n)$ in terms of ${\mathbb P}_{n+1}.$  This, together with  L.~Schl\"afli's famous upper bound, give us asymptotics 
$$P(2, n) \thicksim 2 {2^n-1 \choose n},\quad n\to \infty.$$
Also, the validity of the long-standing conjecture concerning ${\mathbb P}_n$ is proved:
$$\mathbb{P}_n \thicksim (n-1)^22^{1-n}, \quad n\to \infty .$$

\bigskip


\end{abstract}

\let\thefootnote\relax\footnote{This research has been partially supported by RFBR grant 18-01-00398 A}

%\tableofcontents

\section{Introduction.}

\begin{df}\label{df1}
A function $f: \{\pm1\}^n \to \{\pm 1 \}$ is called a threshold function, if there exist real numbers $\alpha_0, \alpha_1, \ldots , \alpha_n,$
such that
$$ f(x_1, \ldots , x_n) =1 \mbox{\quad iff \quad} \alpha_1x_1+ \cdots + \alpha_n x_n + \alpha_0 \geq 0.$$
\end{df}

\noindent Denote by $P(2,n)$ the number of threshold functions.

Let us note that
$$f(x_1,\ldots , x_n) = \mbox{sign}  \langle \bar \alpha , (1, \bar x) \rangle,$$
where \ $(1, \bar x) = (1, x_1, \ldots, x_n) \in {\bf R}^{n+1} \mbox{\quad and \quad} \bar \alpha = (\alpha_0, \ldots , \alpha_n) \in {\bf R}^{n+1}.$
\noindent This observation allows us to correspond a threshold function its $(n+1)$-weight vector $\bar \alpha$ as a point in the dual space $({\bf R}^{n+1})^{\ast} = {\bf  R}^{n+1}$. 

Let $A^{\perp}$ be a finite arrangement of hyperplanes all passing through the zero in ${\bf R}^{n+1}$ (central arrangement) and denote by $A = \{ w_1, \ldots , w_T \}$ the set of their normal vectors. For any $w \in {\bf R}^{n+1}\setminus {\bf 0},$ we consider the linear space $\langle w \rangle,$ generated by $w,$ as a point of the projective space ${\bf RP}^n$. By definition, two hyperplane arrangements  
$A_1 = \{ w_1^1, \ldots , w_T^1 \}$ and  $A_2 = \{ w_1^2, \ldots , w_S^2 \}$ are equal, $A_1^{\perp} \equiv A_2^{\perp}$, iff subsets $\langle A_1 \rangle \MYdef \{\langle w_1^1\rangle, \ldots , \langle w_T^1\rangle \} \subset {\bf RP}^n$ and 
$\langle A_2 \rangle \MYdef \{\langle w_1^2\rangle, \ldots , \langle w_S^2\rangle \} \subset {\bf RP}^n$ coinside, $\langle A_1\rangle = \langle A_2 \rangle$.

It is shown in the paper {\rm \cite{Wi2}}, that $P(2,n)$ can be expressed by the number $C(\langle E_n \rangle)$ of disjoint chambers, obtained as compliment in ${\bf R}^{n+1}$ to the arrangement of $2^n$ hyperplanes all passing through the origin with the normal vectors from the set

\begin{equation}\label{eq1}
E_n=\{ (1, b_1, \ldots , b_n) \mid \; b_i \in \{\pm 1 \} ,\; i=1, \ldots,  n \}.
\end{equation}

The upper bound of the number $C(\langle H \rangle)$ for any central arrangement of hyperplanes with a set $H$ of normal vectors was establisched by L.~Schl\"afli in {\rm \cite{Sch}}. For the case $H=E_n,$ we have the following upper bound:

\begin{equation}\label{eq2}
P(2, n) = C(\langle E_n \rangle) \leq 2 \sum_{i=0}^n {2^n-1 \choose i}.
\end{equation}

It should be noted, that in the early 60s of the $20^{th}$ century the upper bound (\ref{eq2}) was obtained by several authors {\rm \cite{Cam}, \cite{Jos}, 
\cite{Wi1}}. The detailed information of contribution of above mentioned authors can also be found in {\rm \cite{Cov}}.

One of the first lower bound of $P(2, n)$ was established by S.~Muroga in {\rm \cite{Mur}}:

\begin{equation}\label{eq3}
P(2, n) \geq 2^{0.33048 n^2}.
\end{equation}

S. Yajima and T. Ibaraki in {\rm \cite{YaI}} improved the order of the logarithm of the lower bound (\ref{eq3}) upto $n^2/2 $ :

\begin{equation}\label{eq4}
P(2, n) \geq 2^{n(n-1)/2 +8} \mbox{\quad for \quad } n \geq 6.
\end{equation}

Further significant improvements of the bound (\ref{eq4}) were obtained basing on the  paper {\rm \cite{Odl}} of A.~M.~Odlyzko. In the paper {\rm \cite{Zue}},  it was noted that from the papers {\rm \cite{Odl}, \cite{Zas}} follows:

\begin{equation}\label{eq5}
C(E) =P(2, n) \geq 2^{n^2 - 10n^2/\ln n +O(n \ln n)} .
\end{equation}

Taking into account the upper bound (\ref{eq2}) and inequality (\ref{eq5}), it is easy to see that

\begin{equation}\label{eq6}
\lim_{n \to \infty}{\frac{\log_2 P(2, n)}{n^2}} = 1.
\end{equation}

 In the paper {\rm \cite{Ir1}},  it was suggested an original geometric construction that, in combination with the result from the paper {\rm \cite{Odl}}, improved the inequality (\ref{eq5}) to:

\begin{equation}\label{eq7}
P(2, n) \geq 2^{n^2\left(1-\frac{7}{\ln n}\right)}\cdot  P\left(2, \left[\frac{7(n-1)}{\log_2 (n-1)}\right]\right) .
\end{equation}

The generalization of the inequality (\ref{eq7}) for the number of threshold $k$-logic functions was obtained in {\rm \cite{IrK}}. Asymptotics of logarithm of the number of polinomial threshold functions has been recently obtained in {\rm \cite{BV}}.

In parallel to finding the asymptotics of the number of threshold functions, studies were conducted to find the asymptotics of the number of singular $\{\pm 1\}$ (or $\{0, 1\}$) $n\times n$-matrices.

Let $M_n=(a_{ij})$ be a random $n\times n $ $\{\pm 1\}$-matrix, whose entries are independent identically distributed (i.i.d.) Bernoulli random variables:
$$\Pr (a_{ij} = 1) = \Pr (a_{ij} = -1) = \frac{1}{2}. $$ 
Many researchers have devoted considerable attention to the old problem of finding the probability
\begin{equation*}\label{eqPr}
\mathbb{P}_n \MYdef \Pr ( \det M_n = 0) 
\end{equation*}
that a random Bernoulli $n\times n $  $\{\pm 1\}$-matrix $M_n$ is singular. 

In 1963, J.~Koml\' os {\rm \cite{Ko1}} proved P.~Erd\" os' conjecture that the probability that a random Bernoulli $n\times n $  $\{0, 1\}$-matrix is singular approaches $0$ as $n$ tends to infinity. It is also true for random Bernoulli \{$\pm 1$\}-matrices:
\begin{equation}\label{eq8}
\mathbb{P}_n = o_n(1).
\end{equation}
In 1977, J.~Komlos {\rm \cite{Ko2}} improved his result by proving that
\begin{equation}\label{eq9}
\mathbb{P}_n < O\left( \frac{1}{\sqrt n}\right).
\end{equation}
The proof of (\ref{eq9}) is based on the lemma usually referred to as the Littlewood-Offord lemma, which was proved by P. Erd\" os ({\rm \cite{Erd}}).

In 1995, J.~Kahn, J.~Koml\' os, and E.~Szemer\' edi  established in {\rm \cite{KKS}} for the first time an exponential decay of the upper bound of the singularity probability of random Bernoulli matrices:
\begin{equation}\label{eq10}
\mathbb{P}_n  \leq (1- \varepsilon + o_n(1) )^n , \mbox{\; where \; } \varepsilon = 0.001.
\end{equation}

In {\rm \cite{TV1}}, T.~Tao and V.~Vu improved the result (\ref{eq10}) for $\varepsilon = 0.06191$, and then in {\rm \cite{TV2}}, they sharpened their technique to prove (\ref{eq10}) for $\varepsilon = 0.25$:
\begin{equation}\label{eq11}
\mathbb{P}_n \leq \left( \frac{3}{4} + o_n(1)\right)^n .
\end{equation}
 
In 2009, Tao-Vu's result (\ref{eq11}) was further improved  by J.~Bourgain, V.~H.~Vu, and P.~M.~Wood (see {\rm \cite{BVW}}).They proved that

\begin{equation}\label{eq12}
\mathbb{P}_n \leq \left(\frac{\sqrt 2 }{2} + o_n(1)\right)^n.
\end{equation}

In 2018, K.~Tikhomirov finally obtained in {\rm \cite{Tik}} that

\begin{equation}\label{eq13}
\mathbb{P}_n = \left(\frac{1 }{2} + o_n(1)\right)^n.
\end{equation}

In this paper, we prove the validity of the long standing conjecture (see {\rm \cite{Ko2}, \cite{Odl}, \cite{KKS}}) that dominant sources of singularity are the cases when a matrix $M_n$ contains two identical (or opposite) rows or two identical (or opposite) columns.

\bigskip 
\noindent {\bf Theorem 6}
{\it Asymptotics of the probability that a random Bernoulli matrix is singular is $(n-1)^22^{1-n}$:
\begin{equation*}\label{ThA}
\mathbb{P}_n \thicksim (n-1)^22^{1-n}, \qquad n\to \infty .
\end{equation*} }
\bigskip

We also obtain a new lower bound for the number of threshold functions
\begin{equation}\label{eq14}
P(2, n) \geq  2\left[ 1- \frac{n^2}{2^n}\left( 1 + o_n(1) \right)\right]\dbinom{2^n-1}{n}.
\end{equation}

Combaining the lower bound (\ref{eq14}) with the upper bound (\ref{eq2}), we get
\bigskip

\noindent {\bf Theorem 7}
{\it Asymptotics of the number of threshold functions is equal to $2{2^n-1 \choose n}:$
\begin{equation*}\label{eq15}
P(2, n) \thicksim 2 {2^n-1 \choose n}, \qquad n\to \infty .
\end{equation*} }
\bigskip

\section{Function ${\bf \eta}^{\bigstar}$ and its properties.}

As we mentioned in the previous section, any central hyperplane arrangement $H^{\perp}$ with the set of normal vectors $H = \{ w_1, \ldots , w_{\bar{T}} \} \in {\bf R}^{n+1}\setminus {\bf 0}$, we can identify with the subset  $\langle H \rangle \MYdef \{\langle w_1\rangle, \ldots , \langle w_T\rangle \} \subset {\bf RP}^n$  of the $n$-dimensional projective space. 
We define a partially ordered set (poset) $L^H$ in the following way. By definition, any subspace of $\R^{n+1}$ generated by some (possibly empty) subset of $H$ is an element of the poset $L^H$. An element $s \in L^H$ is less than an element $t \in L^H$ iff the subspace $t$ contains the subspace $s$. For any poset $P$, we can define a {\it simplicial complex} $\Delta_P$ in the following way. The set of vertices of $\Delta_P$ coincides with the set of elements $P$ and a set of vertices of $P$ defines a simplex of $\Delta_P$ iff this set forms a chain in $P$. Let us denote by $\Delta_{L^H}$ the simplicial complex of the poset
$$\left( 0_{L^H} , 1_{L^H}\right)  \MYdef \{ z \in L^H \mid 0_{L^H} < z < 1_{L^H} \} ,$$
\noindent where $0_{L^H}$ and $1_{L^H}$ are the elements of the poset $L^H$ corresponding to the zero subspace of $\R^{n+1}$ and the subspace $ span \ \langle w_1, \ldots ,w_T\rangle$, respectively. 
Without loss of generality, we can assume that
\begin{equation*}\label{eq17A}
dim \langle H \rangle \MYdef dim\  span \ \langle w_1, \ldots ,w_T\rangle =n+1,
\end{equation*}
i.e.,
\begin{equation*}\label{eq17B}
span \ \langle w_1, \ldots ,w_T\rangle = {\bf R}^{n+1}.
\end{equation*}
It has been shown in {\rm \cite{Zas}} that the number $C(\langle H \rangle)$ of $(n+1)$-dimensional regions into which $\R^{n+1}$ is divided by hyperplanes from the set $H^{\perp}$ can be found by the formula:
\begin{equation}\label{eq16}
C(\langle H \rangle) = \sum_{t \in L^H}{\left|  \mu(0_{L^H} , t)\right|},
\end{equation}
\noindent where $\mu (s , t)$ is M\"obius function of the poset $L^H$.
M\"obius function of partially ordered set in Zaslavsky's formula (\ref{eq16}) for calculation of the number of chambers $C(\langle H \rangle)$ can be interpreted by tools of algebraic topology in the following way.
First, we introduce a simplicial compex $K^H$. The set of vertices of $K^H$ coincides with the set $\langle H \rangle$. A subset $\{\langle w_{i_1} \rangle , \ldots , \langle w_{i_s} \rangle \}$ of $\langle H \rangle$
forms a simplex of $K^H$ iff
\begin{equation*}\label{eq17}
span \ \langle w_{i_1}, \ldots ,w_{i_s}\rangle \ne span \ \langle w_1, \ldots ,w_T\rangle = {\bf R}^{n+1}.
\end{equation*}
Taking into account the results of the papers {\rm \cite{Fol}, \cite{Hal}}, it is possible to show (see {\rm \cite{Ir3}}) that the absolute value of the M\"obius function $|\mu (0_{L^H}, u)|$  is equal to the dimension of the reduced homology group of the complex $K^{H\cap u}$ with coefficients in an arbitrary field $\F$:
\begin{equation}\label{eq18}
|\mu (0_{L^H}, u)| = \rank \, \tilde H_{dim\, u -2}\left( K^{ H  \cap u} ; \F\right).
\end{equation}
\noindent Here, the set $\langle H \rangle \cap u$ consists of all elements $\langle H \rangle$ belonging to the subspace $u \subset \R^{n+1}$ and is considered as a subset of $\R^{dim\, u}\MYdef u.$

It follows from the definition of M\"obius function that
\begin{equation*}\label{eq19}
\sum_{0_{L^H} \leq u < 1_{L^H}}{|\mu (0_{L^H} , u)|} \geq \left| -\sum_{0_{L^H} \leq u < 1_{L^H}}{\mu (0_{L^H} , u)}\right| = |\mu (0_{L^H} , 1_{L^H})|.
\end{equation*}
\noindent Hence,

\begin{equation}\label{eq20}
C(\langle H \rangle) = |\mu (0_{L^H} , 1_{L^H})| + \sum_{0_{L^H} \leq u < 1_{L^H}}{|\mu (0_{L^H} , u)|} \geq 2|\mu (0_{L^H}, 1_{L^H})|.
\end{equation}

\noindent From (\ref{eq18}) and (\ref{eq20}), we have:
\begin{equation}\label{eq21}
C(\langle H \rangle) \geq 2\,\rank \, H_{n-1}\left( K^H ; \F\right).
\end{equation}

\noindent As a consequence of (\ref{eq21}) for the case $H=E_n$, we have:
\begin{equation}\label{eq22}
P(2 , n) = C(\langle E_n \rangle) \geq 2\, \rank\, H_{n-1}\left( K^{E_n} ; \F\right).
\end{equation}

Let us fix an arbitrary order on the set $\langle H \rangle :$
\begin{equation}\label{eq23}
\pi : [T] \to \langle H \rangle \subset {\bf RP}^n, \quad |\langle H \rangle | = T, \quad  \langle w_i \rangle \MYdef \pi(i), \quad 1 \leq i \leq T.
\end{equation}
 
Let us denote by $\langle H \rangle^{\times s}$, $s=1,\ldots, T$, the set of ordered collections $(\langle w_{i_1}\rangle , \ldots ,\langle w_{i_s} \rangle)$ of different $s$  elements from $\langle H \rangle$ and let $\langle H \rangle^{\times s}_{\ne 0}\subset \langle H \rangle^{\times s}$ and $\langle H \rangle^{\times s}_{= 0}\subset \langle H \rangle^{\times s}$ be the subsets 
\begin{equation}\label{eq23A}
\langle H\rangle ^{\times s}_{\ne 0} \MYdef \{(\langle w_{i_1}\rangle, \ldots ,\langle w_{i_s}\rangle) \in \langle H \rangle^{\times s} \mid dim \, span \ \langle w_{i_1}, \ldots , w_{i_s} \rangle  = s \}.
\end{equation}
\begin{equation}\label{eq23B}
\langle H\rangle ^{\times s}_{= 0} \MYdef \{(\langle w_{i_1}\rangle, \ldots ,\langle w_{i_s}\rangle) \in \langle H \rangle^{\times s} \mid dim \, span \ \langle w_{i_1}, \ldots , w_{i_s} \rangle  < s \}.
\end{equation}

\begin{df}\label{df2}
We say that an ordered collection of different elements $(\langle w_{i_1} \rangle, \ldots , \langle w_{i_n} \rangle ) \in \langle H \rangle^{\times n}$  satisfies to $\eta^{\pi}_n(\langle H \rangle)$ condition iff the following requirements are fullfilled: 
\begin{enumerate}
	\item $2\leq i_1 < i_2 < \cdots < i_n \leq T;$
	\item $\forall l , \quad 1\leq l\leq n,$ the element $\langle w_{i_l} \rangle$ is minimal in  order $\pi$  among all points from the set $\langle H \rangle \bigcap span \ \langle \langle w_{i_l}\rangle, \ldots , \langle w_{i_n}\rangle \rangle.$
\end{enumerate}
\end{df}

It follows from the definition \ref{df2} that if a collection $(\langle w_{i_1} \rangle, \ldots , \langle w_{i_n} \rangle )$ satisfies to $\eta^{\pi}_n(\langle H \rangle)$ condition, then for all $l=1, \ldots, n,$ we have 
\begin{equation*}\label{eq24}
dim \ span \ \langle w_{i_l}, \ldots, w_{i_n} \rangle = n- l+1, 
\end{equation*}
i.e.,
\begin{equation*}\label{eq24A}
(\langle w_{i_1} \rangle, \ldots , \langle w_{i_n} \rangle ) \in \langle H\rangle ^{\times n}_{\ne 0},
\end{equation*}
and
\begin{equation}\label{eq25}
\ span \ \langle w_1,  w_{i_1}, \ldots, w_{i_n} \rangle = {\bf R}^{n+1}. 
\end{equation}
Denote by $B^{\pi}(\langle H \rangle)$ the set
\begin{equation}\label{eq25A}
B^{\pi}(\langle H \rangle)\MYdef \{W \in  \langle H\rangle ^{\times n}_{\ne 0} \mid W \mbox{satisfyes to} \  \eta^{\pi}_n(\langle H \rangle)\  \mbox{condition}\}.
\end{equation}

The theorem~7 of {\rm \cite{Ir3}} is also true for any finite subset $\langle H \rangle \subset {\bf RP}^n$. It asserts that the number of collections $(\langle w_{i_1} \rangle, \ldots , \langle w_{i_n} \rangle )$ satisfying to $\eta^{\pi}_n(\langle H \rangle)$ condition is equal to the $ rank\ H_{n-1}(K^H;{\bf F})$ . Hence, the number of collection satisfying to $\eta^{\pi}_n(\langle H \rangle)$ condition doesn't depend on the order $\pi$ on the set $\langle H \rangle.$ Let us denote this number by $\eta^{\bigstar}_n(\langle H \rangle).$
Thus on the set $2^{{\bf RP}^n}_{fin}$ of finite subsets of ${\bf RP}^n,$ the function $\eta^{\bigstar}_n : 2^{{\bf RP}^n}_{fin} \to \mathbb{Z}_{\geq 0}$ satisfies to the formula:
\begin{equation}\label{eq26}
\eta^{\bigstar}_n(\langle H \rangle) = rank\ H_{n-1}(K^H;{\bf F}), \quad \langle H \rangle \subset {\bf RP}^n.
\end{equation}

\begin{pr}\label{pr1}
$\eta^{\bigstar}_n$ is a supermodular function on $2^{{\bf RP}^n}_{fin}$.
\end{pr}
\proof It is necessary to demonstrate that for any finite subset $\langle H \rangle \subset {\bf RP}^n$, $|\langle H \rangle| =T,$ and any two different elements $\langle u\rangle, \langle v \rangle \in {\bf RP}^n \setminus \langle H \rangle $ the following inequality 
\begin{equation}\label{eq27}
\eta^{\bigstar}_n(\langle H \rangle \cup \{\langle u \rangle \}) - \eta^{\bigstar}_n(\langle H \rangle) \leq \eta^{\bigstar}_n(\langle H \rangle \cup \{\langle u \rangle, \langle v \rangle \}) - \eta^{\bigstar}_n(\langle H \rangle \cup \{\langle v \rangle \}) 
\end{equation}
is fullfiled.

For any order $\pi : [T] \to \langle H \rangle$, we define orders $\pi^{u,v} : [T+2] \to  \langle H \rangle \cup \{\langle u \rangle, \langle v \rangle \}$, $\pi^{u} : [T+1] \to  \langle H \rangle \cup \{\langle u \rangle \},$ and 
$\pi^{v} : [T+1] \to  \langle H \rangle \cup \{\langle v \rangle \}$ such that 
\begin{equation}\label{eq28}
\begin{split}
\pi^{u,v}(i) =\pi^u(i) = \pi^v(i) = \pi (i), \quad \forall i=1, \ldots , T;\\
\pi^{u,v} (T+1) = \langle u \rangle, \quad \pi^{u,v}(T+2) = \langle v \rangle ; \\
\pi^{u} (T+1) = \langle u \rangle, \quad \pi^{v} (T+1) = \langle v \rangle.
\end{split}
\end{equation}
Then the expression in the left part of the inequality (\ref{eq27}) equals to the number of collections 
$(\langle w_{i_1} \rangle, \ldots, \langle w_{i_{n-1}} \rangle, \langle u \rangle)$ satisfying to $\eta^{\pi^u}_n(\langle H \rangle \cup \{\langle u \rangle \})$ condition. Due to (\ref{eq28}), these collections also satisfy to 
$\eta^{\pi^{u, v}}_n(\langle H \rangle \cup \{\langle u \rangle, \langle v \rangle \})$ condition. The expression in the right side of the inequality (\ref{eq27}) equals to cardinality of the set consisting of collections of the form 
$(\langle w_{i_1} \rangle, \ldots, \langle w_{i_{n-1}} \rangle, \langle u \rangle)$ and $(\langle w_{j_1} \rangle, \ldots, \langle w_{j_{n-2}} \rangle, \langle u \rangle, \langle v \rangle)$ satisfying to 
$\eta^{\pi^{u, v}}_n(\langle H \rangle \cup \{\langle u \rangle, \langle v \rangle \})$ condition. Hence, the inequality (\ref{eq27}) is proved.

\begin{flushright} {\sc Q.E.D.} \end{flushright}

\bigskip

Denote by $P^{\langle w \rangle}_{R^{n+1}}$ the orthogonal projector along the linear subspace $\langle w \rangle \subset {\bf R}^{n+1}$ onto its $n$-dimensional orthogonal compliment $\langle w \rangle^{\perp} \subset {\bf R}^{n+1},$ and denote by $v^{\perp w}$ the image of a vector $v \in {\bf R}^{n+1}:$
\begin{equation}\label{eq29}
v^{\perp w} \MYdef P^{\langle w \rangle}_{R^{n+1}}(v).
\end{equation}

For $\langle H \rangle = \{\langle w_1\rangle, \ldots , \langle w_T\rangle \} \subset {\bf RP}^n$ and  $\langle w \rangle \not\in \langle H \rangle,$ we denote by $\langle H \rangle^{\perp w}$ the set:
\begin{equation}\label{eq30}
\langle H \rangle^{\perp w} \MYdef  P^{\langle w \rangle}_{R^{n+1}}(\langle H \rangle) = \{\langle w_1^{\perp w} \rangle, \ldots , \langle w_T^{\perp w} \rangle \} \subset {\bf RP}^{n-1}.
\end{equation}

\begin{thm}\label{thm1}
For any finite subset $\langle H \rangle \subset {\bf RP}^n$ and element $\langle u \rangle \in {\bf RP}^n \setminus \langle H \rangle,$ we have:
\begin{equation}\label{eq31}
\eta^{\bigstar}_n(\langle H \rangle \cup \{\langle u \rangle \})  =  \eta^{\bigstar}_n(\langle H \rangle) + \eta^{\bigstar}_{n-1}(\langle H \rangle^{\perp u}) .
\end{equation}
\end{thm}

\proof Let $\pi : [T+1] \to  \langle H \rangle \cup \{\langle u \rangle \}$ be an order on $\langle H \rangle \cup \{\langle u \rangle \} \subset {\bf RP}^n$ such that $\pi(T+1) = \langle u \rangle$.
For any $\langle w^{\perp u} \rangle \in \langle H \rangle^{\perp u},$ let
\begin{equation}\label{eq31A}
\begin{split}
T^{\pi}(w^{\perp u}) \MYdef \{ i\in [T] \mid \langle w_i^{\perp u}\rangle = \langle w^{\perp u} \rangle\} \\
m(w^{\perp u}) \MYdef min \ \{i\in T^{\pi}(w^{\perp u})\}
\end{split}
\end{equation}
For $\langle x^{\perp u} \rangle, \langle y^{\perp u} \rangle \in \langle H\rangle^{\perp u},$ we say that 
\begin{equation}\label{eq31B}
\langle x^{\perp u} \rangle <_{\pi^{u^{\perp}}} \langle y^{\perp u}\rangle \quad \mbox{iff} \quad m(x^{\perp u}) < m(y^{\perp u}).
\end{equation}
Let $|\langle H\rangle^{\perp u}| =T'.$ Then we define the order $\pi|_{u^{\perp}} : [T'] \to \langle H\rangle^{\perp u}$ as the unique map preserving the linear orders:
\begin{equation}\label{eq32}
i<j \quad \mbox{iff} \quad \pi|_{u^{\perp}}(i) <_{\pi^{u^{\perp}}} \pi|_{u^{\perp}}(j).
\end{equation}

In the proof of the Proposition \ref{pr1}, we have shown that the cardinality of the set  
\begin{multline*}
B^{\pi}_u(\langle H \rangle \cup \{\langle u \rangle \})\MYdef \\
=\left\{W \in  B^{\pi}(\langle H \rangle \cup \{\langle u \rangle \}) \mid W = (W', \langle u \rangle ), \ W' \in \langle H\rangle ^{\times (n-1)}_{\ne 0}\right\}
\end{multline*}
\noindent is equal to the number $\eta^{\bigstar}_n(\langle H \rangle \cup \{\langle u \rangle \}) - \eta^{\bigstar}_n(\langle H \rangle):$
\begin{equation}\label{eq33}
\eta^{\bigstar}_n(\langle H \rangle \cup \{\langle u \rangle \}) - \eta^{\bigstar}_n(\langle H \rangle) =|B^{\pi}_u(\langle H \rangle \cup \{\langle u \rangle \})|.
\end{equation}

For any $W=(\langle w_{i_1} \rangle, \ldots , \langle w_{i_{n-1}} \rangle, \langle u \rangle) \in B^{\pi}_u(\langle H \rangle \cup \{\langle u \rangle\}),$ we assert that
\begin{equation}\label{eq34}
W^{\perp u} \MYdef (\langle w_{i_1}^{\perp u} \rangle, \ldots , \langle w_{i_{n-1}}^{\perp u} \rangle ) \in B^{\pi|_{u^{\perp}}}\left(\langle H\rangle^{\perp u}\right).
\end{equation}
First of all, we note that for $W\in  B^{\pi}_u(\langle H \rangle \cup \{\langle u \rangle\}),$ we have 
\begin{equation}\label{eq34A}
i_l =  m(i_l) \MYdef m(w_{i_l}^{\perp u}), \quad \forall l=1, \ldots, n-1.
\end{equation}
Indeed, if the condition (\ref{eq34A}) is not fulfilled for some $l, 1\leq l \leq n-1$, then $i_l > m(i_l)$. 
Taking into account that 
\begin{equation}\label{eq34B}
w_i^{\perp u} = w_i - \beta_i u, \quad \beta_i=\frac{(w_i, u)}{(u, u)}, \quad i=1, \ldots, n-1,
\end{equation}
we have 
\begin{equation}\label{eq34C}
\langle w_{m(i_l)} \rangle \in span \ \langle \langle w_{i_l} \rangle, \langle u\rangle \rangle \subset  span \ \langle \langle w_{i_l} \rangle, \ldots, \langle w_{n-1} \rangle ,\langle u\rangle \rangle .
\end{equation}
The inclusion (\ref{eq34C}) contradicts to our choice $W\in  B^{\pi}_u(\langle H \rangle \cup \{\langle u \rangle\}).$

Now consider the case when the condition (\ref{eq34}) is not fulfilled, i.e., there exist $l, k$, $1\leq l \leq n-1,$ and $k<i_l$ such that
\begin{equation*}\label{eq35}
\langle w_k^{\perp u}\rangle \in span \ \langle \langle w_{i_l}^{\perp u} \rangle. \ldots , \langle w_{i_{n-1}}^{\perp u} \rangle \rangle.
\end{equation*}
Hence, there exist $\alpha_l, \ldots, \alpha_{n-1} \in {\bf R}$ such that
\begin{equation}\label{eq36}
w_k^{\perp u} =\alpha_l w_{i_l}^{\perp u} + \cdots + \alpha_{n-1} w_{i_{n-1}}^{\perp u}.
\end{equation}
From (\ref{eq34B}) and (\ref{eq36}) we have
\begin{equation}\label{eq37}
w_k =\alpha_l w_{i_l} + \cdots + \alpha_{n-1} w_{i_{n-1}} - (\alpha_l \beta_{i_l} + \cdots + \alpha_{n-1} \beta_{i_{n-1}} - \beta_k)u.
\end{equation}
The last equation means that $\langle w_k\rangle \in span \ \langle \langle w_{i_l} \rangle, \ldots, \langle w_{i_{n-1}} \rangle, \langle u \rangle \rangle$, i.e., conradicts to our requirement $W\in  B^{\pi}_u(\langle H \rangle \cup \{\langle u \rangle\}).$ 

Thus we can define the map 
$$\psi^{\pi}_u:  B^{\pi}_u(\langle H \rangle \cup \{\langle u \rangle\}) \to  B^{\pi|_{u^{\perp}}}\left(\langle H\rangle^{\perp u}\right)$$
by the rule
\begin{equation}\label{eq38}
\psi^{\pi}_u (W) \MYdef W^{\perp u}
\end{equation}
We assert that $\psi^{\pi}_u$ is injective. If we assume the opposite, then there exist $W_1, W_2 \in  B^{\pi}_u\left(\langle H \rangle \cup \{\langle u \rangle\right\}),$ $W_1\ne W_2$, $W_1=(\langle w_{i_1} \rangle, \ldots , \langle w_{i_{n-1}}\rangle, \langle u \rangle),$
$W_2=(\langle w_{j_1} \rangle, \ldots , \langle w_{j_{n-1}}\rangle, \langle u \rangle)$ such that
\begin{equation}\label{eq39}
\psi^{\pi}_u (W_1)=(\langle w_{i_1}^{\perp u} \rangle, \ldots , \langle w_{i_{n-1}}^{\perp u} \rangle ) = (\langle w_{j_1}^{\perp u} \rangle, \ldots , \langle w_{j_{n-1}}^{\perp u} \rangle ) = \psi^{\pi}_u (W_2).
\end{equation}
Let $l$, $1\leq l\leq n-1,$ be the maximal number such that $w_{i_l}\ne w_{j_l}.$ Without loss of generality, we can assume that $i_l < j_l.$ It follows from (\ref{eq34A}) that $m(i_l)=i_l<j_l=m(j_l).$ Then 
\begin{equation}\label{eq40}
\langle w_{i_l}^{\perp u} \rangle = \langle w_{m(i_l)}^{\perp u} \rangle <_{\pi^{u^{\perp}}} \langle w_{m(j_l)}^{\perp u} \rangle = \langle w_{j_l}^{\perp u} \rangle.
\end{equation}
The inequality (\ref{eq40}) contradicts to our assumption (\ref{eq39}) that $\langle w_{i_l}^{\perp u} \rangle=\langle w_{j_l}^{\perp u} \rangle,$ $l=1, \ldots, n-1.$

We assert that $\psi^{\pi}_u$ is surjective. 
Let $ X^{\perp u} =(\langle x_1 \rangle, \ldots , \langle x_{n-1} \rangle ) \in B^{\pi|_{u^{\perp}}}\left(\langle H\rangle^{\perp u}\right).$ Let us put
\begin{equation}\label{eq41}
\left(\psi^{\pi}_u\right)^{-1}\left(X^{\perp}\right) \MYdef (\langle w_{m(x_1)} \rangle, \ldots, \langle w_{m(x_{n-1})} \rangle, \langle u \rangle).
\end{equation}
It is necessary to demonstrate that 
\begin{equation}\label{eq42}
\left(\psi^{\pi}_u\right)^{-1}\left(X^{\perp}\right)  \in B^{\pi}_u(\langle H \rangle \cup \{\langle u \rangle \}).
\end{equation}
If we assume that the inclusion (\ref{eq42}) isn't true, then there exist $l, k,$ $1\leq l \leq n-1,$ and $k< m(x_l)$ such that
\begin{equation*}\label{eq43}
\langle w_k\rangle \in span \ (\langle w_{m(x_l)}\rangle, \ldots , \langle w_{m(x_{n-1})} \rangle, \langle u \rangle ).
\end{equation*}
Hence,
\begin{equation}\label{eq44}
\langle w_k^{\perp u}\rangle \in span \ (\langle w_{m(x_l)}^{\perp u}\rangle, \ldots , \langle w_{m(x_{n-1})}^{\perp u} \rangle ).
\end{equation}
Since $k<m(x_l)$, then 
\begin{equation}\label{eq45}
m(w_k^{\perp u}) < m(x_l) < \cdots < m(x_{n-1}).
\end{equation}
From inclusion (\ref{eq44}) and inequalities (\ref{eq45}) we get a contradiction to  $ X^{\perp u} \in B^{\pi|_{u^{\perp}}}\left(\langle H\rangle^{\perp u}\right).$

Thus we have demonstrated that
\begin{equation}\label{eq46}
|B^{\pi}_u(\langle H \rangle \cup \{\langle u \rangle \})| = |B^{\pi|_{u^{\perp}}}\left(\langle H\rangle^{\perp u}\right)|.
\end{equation}

Since
\begin{equation}\label{eq47}
\eta^{\bigstar}_{n-1} (\langle H\rangle^{\perp u}) = |B^{\pi|_{u^{\perp}}}\left(\langle H\rangle^{\perp u}\right)|,
\end{equation}
\noindent our Theorem follows from the equalities (\ref{eq33}),  (\ref{eq46}), and  (\ref{eq47}).

\begin{flushright} {\sc Q.E.D.} \end{flushright}

For any finite subset $\langle H \rangle \subset {\bf RP}^n$ and element $\langle w \rangle \in {\bf RP}^n ,$ we denote by $\dbinom{\langle H \rangle }{ \eta^{\bigstar}_n}^{\langle w \rangle}$ the following sum:
\begin{equation}\label{eq48}
\dbinom{\langle H \rangle }{ \eta^{\bigstar}_n}^{\langle w \rangle} \MYdef \sum_{\{\langle w_{i_1} \rangle, \ldots, \langle w_{i_n} \rangle\}\subset \langle H \rangle} \eta^{\bigstar}_n \left(\{\langle w \rangle, \langle w_{i_1}\rangle, \ldots, \langle w_{i_n} \rangle\}\right).
\end{equation}

\begin{thm}\label{thm2}
For any $n\geq 1$, finite subset $\langle H \rangle \subset {\bf RP}^n,$ and element $\langle w \rangle \in {\bf RP}^n,$ we have:
\begin{equation}\label{eq49}
\eta^{\bigstar}_{n} (\langle H\rangle \cup \{\langle w \rangle \}) \leq \dbinom{\langle H \rangle }{ \eta^{\bigstar}_n}^{\langle w \rangle}.
\end{equation}
\end{thm}

\proof  Let $\pi : [T] \to  \langle H \rangle \cup \{\langle w \rangle \}$ be an order on $\langle H \rangle \cup \{\langle w \rangle \} \subset {\bf RP}^n$, $T=|\langle H \rangle \cup \{\langle w \rangle \}|,$  such that $\pi(1) = \langle w \rangle$. It is easy to see that 
\begin{equation*}\label{eq52}
\begin{split}
\eta^{\bigstar}_n \left(\{\langle w \rangle, \langle w_{i_1} \rangle, \ldots, \langle w_{i_n} \rangle\}\right)=1 \Leftrightarrow \\
\Leftrightarrow span \ \langle w, w_{i_1}, \ldots, w_{i_n}\rangle = {\bf R}^{n+1}.
\end{split}
\end{equation*}
It follows from definition \ref{df2} that if a collection $(\langle w_{i_1} \rangle, \ldots , \langle w_{i_n} \rangle ) \in B^{\pi}(\langle H \rangle \cup \{\langle w \rangle \}),$ then (see (\ref{eq25})) 

$$\ span \ \langle w_1,  w_{i_1}, \ldots, w_{i_n} \rangle =  span \ \langle w, w_{i_1}, \ldots, w_{i_n}\rangle = {\bf R}^{n+1}.$$ 

Now the Theorem follows from the equality $\eta^{\bigstar}_n(\langle H \rangle \cup \{\langle w \rangle \}) = | B^{\pi}(\langle H \rangle \cup \{\langle w \rangle \})|.$ 

\begin{flushright} {\sc Q.E.D.} \end{flushright}

\bigskip

\section{A formula for ${\bf \eta}^{\bigstar}_*$ in terms of combinatorial flags on a central hyperplane arrangement.}

For any $W= (\langle w_{i_1} \rangle, \ldots , \langle w_{i_n}\rangle) \in \langle H \rangle^{\times n}_{\ne 0}$  and $ l = 1, \ldots, n,$ let 
\begin{equation}\label{eq85}
q^W_l \MYdef |L_l(W) \cap \langle H \rangle| \MYdef |span \ \langle \langle w_{i_{n-l+1}} \rangle, \ldots , \langle w_{i_n} \rangle \rangle \cap \langle H \rangle|.
\end{equation}

\begin{df}\label{df3}
For any $W \in \langle H\rangle^{\times n},$ the ordered set of numbers
\begin{equation}\label{eq86}
W(\langle H \rangle) \MYdef (q^W_n, q^W_{n-1}, \ldots, q^W_1)
\end{equation}
is called a combinatorial flag on $\langle H \rangle \subset {\bf RP}^n$ of the ordered set $W$.

If $W \in \langle H \rangle^{\times n}_{\ne 0}$, then $W(\langle H \rangle)$ is called a full combinatorial flag of $W$.
\end{df}
For the sake of simplicity, we will use the following notation:
\begin{equation}\label{eq87}
W[H] \MYdef q^W_n\cdot q^W_{n-1} \cdots q^W_1.
\end{equation}

To define the set $B^{\pi}(\langle H \rangle)$ (see (\ref{eq25A})), we fixed an order $\pi : [T] \to \langle H \rangle \subset {\bf RP}^n$ (see (\ref{eq23})) that allowed us to compare elements of $\langle H \rangle:$
$$  \langle w_i \rangle <_{\pi} \langle  w_j \rangle \iff i<j.$$

We denote by $\Gamma$ the set of all orders on the set $\langle H \rangle$. Then any order on $\langle H \rangle$ can be defined as composition
$$[T] \mapr{\gamma} [T] \mapr {\pi} H$$
\noindent of a permutation $\gamma: [T] \to [T]$ with $\pi$, and 
\begin{multline}\label{eq88}
\langle w_i \rangle <_{\gamma} \langle w_j \rangle \iff (\pi \gamma)^{-1}(\langle w_i \rangle) < (\pi \gamma)^{-1}(\langle w_j \rangle) \iff   \\ 
\iff \gamma^{-1}(i) < \gamma^{-1}(j).
\end{multline}
Thus $\Gamma$ can be identified with the symmetric group $Sym([T]),$ and any permutation $\sigma : [T] \to [T]$ defines the basis of the homology group $H_{n-1} \left(K^H ; {\bf F}\right),$ considered as a vector space over ane fixed field ${\bf F},$ say ${\mathbb Z}_2,$ as the subset of collections of $n$ elements from $\langle H \rangle$
$$B^{\pi \circ \sigma}(\langle H \rangle) \subset \langle H\rangle^{\times s}_{\ne 0} $$
\noindent obeying to $\eta^{\pi\circ \sigma}_n(\langle H \rangle)$ condition.

\begin{thm}\label{thm4}
For any probability distribution $p=(p_1, \ldots , p_T)$ on a subset $\langle H \rangle \subset {\bf RP}^n,$ $span \ \langle H \rangle = {\bf R}^{n+1},$ the following equality is true: 
\begin{equation}\label{eq89}
\eta^{\bigstar}_n(\langle H \rangle ) = \sum_{W\in \langle H\rangle^{\times n}_{\ne 0}}{\frac{1- p_{i_1} -p_{i_2} - \cdots - p_{i_{q_n^W}}}{W[H]}}.
\end{equation}
Here, the indices used in the numerator correspond to elements from 
$$L_n(W)\cap \langle H \rangle = \left\{\langle w_{i_1} \rangle , \ldots, \langle w_{i_n} \rangle, \ldots , \langle w_{i_{q_n^W}}\rangle \right\}.$$
\end{thm}

\proof  We define the probability distribution $\tilde p$ on the set $\Gamma \cong Sym([T])$ by the formula:
\begin{equation*}\label{eq90}
\tilde p(\gamma) =p_{\gamma(1)}\frac{1}{(T-1)!} , \quad \gamma \in Sym([T]).
\end{equation*}

For any collection $W=(\langle w_{i_1} \rangle, \ldots , \langle w_{i_n}) \rangle \in H^{\times n}_{\ne 0},$ we define the random function $I_W : \Gamma \to \R$ by the formula:
\begin{equation*}\label{eq91}
I_W(\gamma) \MYdef \left\{\begin{array}{ll}
1, & \mbox{if}\  W \  \mbox{satisfies to} \ \eta^{\pi\circ \gamma}_n(\langle H \rangle)\ \mbox{condition}; \\  
    & \\
0, & \mbox{in all other cases.}
\end{array}\right. 
\end{equation*}

Let
\begin{equation*}\label{eq92}
I \MYdef \sum_{W\in \langle H\rangle^{\times n}_{\ne 0}}{I_W} : \Gamma \to \R.
\end{equation*}

Then for any $\gamma \in \Gamma,$
\begin{equation*}\label{eq93}
I(\gamma) = \const = |B^{\pi \circ \gamma}(\langle H \rangle)| = \rank \, H_{n-1}\left(K^H ; \F\right) = \eta^{\bigstar}_n(\langle H \rangle ).
\end{equation*}
Hence, the expectation of $I$ is equal to $\eta^{\bigstar}_n(\langle H \rangle )$ :
\begin{equation}\label{eq94}
\mathbb{E} [I] = \eta^{\bigstar}_n(\langle H \rangle ) .
\end{equation}
Additivity of expectation reduces the problem of calculation $\mathbb{E}[I]$ to counting the probability $\Pr (I_W =1)$:
\begin{equation}\label{eq95}
\mathbb{E}[I] = \sum_{W\in \langle H\rangle^{\times n}_{\ne 0}}{\mathbb{E}[I_W]} = \sum_{W\in \langle H\rangle^{\times n}_{\ne 0}}{\Pr (I_W =1)}.
\end{equation}

Further we calculate the number of permutations $\gamma$ such that $I_W(\gamma)=1$. Since $\langle w_{\gamma(1)}\rangle \notin L_n(W),$ then $q^W_n$ elements from 
$L_n(W)\cap \langle H \rangle$ can be located in any places except the first one, i.e., $\gamma^{-1}(j) \ne 1$ for any $j\in [T]$ such that $\langle w_j \rangle  \in L_n(W)\cap \langle H \rangle.$ The arrangement of the remaining elements from 
$\langle H \rangle \setminus \{\langle w_{\gamma (1)}\rangle \cup \{L_n(W)\cap \langle H \rangle \}\}$ does not affect the fulfillment of $\eta^{\pi\circ \gamma}_n(\langle H \rangle)$  condition.  By $\eta^{\pi\circ \gamma}_n(\langle H \rangle)$ condition, the element $\langle w_{i_1} \rangle$ has to be in the first place among the selected $q^W_n$ positions for arrangement of the set $L_n(W)\cap \langle H\rangle,$ while $q_{n-1}^W$ elements from $L_{n-1}(W)\cap \langle H\rangle$ can be located in any of the remained $q_n^W-1$ places. The arrangement of the elements from $\{L_n(W)\cap \langle H \rangle\} \setminus \{\langle w_{i_1}\rangle \cup \{L_{n-1}(W)\cap \langle H \rangle\}\}$ in $q_n^W-q_{n-1}^W-1$ places, left after choosing $q_{n-1}^W+1$ places for arrangement of the set $L_{n-1}(W)\cap \langle H \rangle$ and $\langle w_{i_1} \rangle,$ doesn't affect the fulfillment of $\eta^{\pi\circ \gamma}_n(\langle H \rangle)$ condition. Continuing the same way, we get that in the first place among $q^W_l$ positions selected for the elements from $L_l(W)\cap \langle H \rangle$ has to be located the element $\langle w_{i_{n-l+1}} \rangle$, while $q^W_{l-1}$ elements from $L_{l-1}(W)\cap \langle H \rangle$ can be located in any of the remained $q^W_l-1$ places, and the positions of the elements from $\{L_l(W)\cap \langle H\rangle\} \setminus \{\langle w_{i_{n-l+1}}\rangle \cup \{L_{l-1}(W)\cap \langle H\rangle \}\}$ in $q_l^W-q_{l-1}^W-1$  places, left after choosing $q_{l-1}^W+1$ places for arrangement of the set $L_{l-1}(W)\cap \langle H \rangle$ and $w_{i_{n-l+1}},$ doesn't affect the fulfillment of $\eta^{\pi\circ \gamma}_n(\langle H \rangle)$ condition.

Denote by $N(\gamma(1)=i)$ the number of permutations $\gamma$ with fixed value $\gamma(1)=i$ such that $\langle w_i \rangle \notin L_n(W)$. Then
$$\begin{array}{l} 
N(\gamma(1)=i) = 
{T-1 \choose q_n^W} (T-1-q_n^W)!\cdot {q_n^W-1 \choose q_{n-1}^W}(q_n^W-q_{n-1}^W-1)!\cdots \\
\\
\cdot {q_l^W-1 \choose q_{l-1}^W}(q_l^W-q_{l-1}^W-1)!\cdots {q_2^W-1 \choose q_1^W}(q_2^W-q_1^W-1)!=  \\
\\
= \frac{(T-1)!}{q_n^W!}\cdot \frac{(q_n^W-1)!}{q_{n-1}^W!}\cdots \frac{(q_l^W-1)!}{q_{l-1}^W!}\cdots \frac{(q_2^W-1)!}{q_1^W!}  
= \frac{(T-1)!}{q_n^W q_{n-1}^W \cdots q_2^W\cdot q_1^W} = \\
\\
=\frac{(T-1)!}{W[H]} \quad (\mbox{since} \ q_1^W =1).
\end{array}$$
\noindent Then we have
$$\Pr (I_W=1) = \sum_{i\in [T]\, s.t. \, \langle w_i \rangle \notin L_n(W)}{p_i\frac{1}{(T-1)!}\frac{(T-1)!}{W[H]}} = $$
$$= \frac{1-p_{i_1} -\cdots -  p_{i_{q_n^W}}}{W[H]},$$
\noindent where $ L_n(W)\cap \langle H \rangle = \{\langle w_{i_1}\rangle ,\ldots, \langle w_{i_n} \rangle, \ldots, \langle w_{i_{q_n^W}}\rangle\}.$
Now the Theorem follows from  (\ref{eq94}) and  (\ref{eq95}).

\begin{flushright} {\sc Q.E.D.} \end{flushright}

\begin{rem}
Since the right side of equation of Theorem \ref{thm4} is expressed by a polynomial of degree 1, then the Theorem \ref{thm4} is true for any $p_i \in {\bf R}$, $i=1,\ldots, T$, such that $ \sum_{i=1}^T{p_i}=1$.
\end{rem}

\vskip 1cm

\section{A lower bound for ${\bf \eta}^{\bigstar}_*$.}

Next, we are going to get an upper bound for the $\rank \, H_n\left( K^H , K^H_{n-1}; \F\right),$ where $\langle H \rangle = \{\langle w_1\rangle, \ldots , \langle w_T\rangle \} \subset {\bf RP}^n$, and $K^H_{n-1}$ is the $(n-1)$-skeleton of $K^H.$
The nonzero part of the homology exact sequence of the pair $(K^H, K^H_{n-1})$ has the following form (see (11) of {\rm \cite{Ir3}}):

\begin{equation}\label{eqC1}
0 \to H_n\left( K^H , K^H_{n-1}; \F\right) \to H_{n-1}\left( K^H_{n-1}; \F\right) \to H_{n-1}\left( K^H ; \F\right) \to 0 .
\end{equation} 

For any $\Delta = (\langle w_{i_1}\rangle, \ldots ,\langle w_{i_{n+1}}\rangle) \in \langle H \rangle^{\times (n+1)}_{=0},$ $i_1< \ldots < i_{n+1},$ let put

\begin{equation}\label{eqC2}
n(\Delta ) \MYdef \{t \in [n] \mid w_{i_t} \in \, span \ \langle w_{i_{t+1}}, \ldots , w_{i_{n+1}} \rangle  \},
\end{equation}

\begin{equation}\label{eqC3}
\begin{split}
\Delta (H) \MYdef \{ \langle w_p \rangle  \in \langle H \rangle \mid \exists  t \in [n] \,  \mbox{s.t.} \, p< i_{t+1}, \mbox{and}\\
w_p \in \, span \ \langle w_{i_{t+1}}, \ldots , w_{i_{n+1}} \rangle  \},
\end{split}
\end{equation}

\begin{equation}\label{eqC4}
t(\Delta ) \MYdef \max_{t \in n(\Delta)} t,
\end{equation}

\begin{equation}\label{eqC5}
\langle w(\Delta) \rangle \MYdef \max_{\langle w \rangle \in \Delta (H)} \langle w \rangle.
\end{equation}

\noindent  Let
\begin{equation}\label{eqC6}
\begin{split}
C^{\pi}_n(H) \MYdef \{\Delta = (\langle w_{i_1}\rangle, \ldots ,\langle w_{i_{n+1}}\rangle) \in \langle H \rangle^{\times (n+1)}_{=0} \mid \\
i_1< \ldots < i_{n+1},  \langle w_{i_{t(\Delta)}} \rangle =\langle w(\Delta) \rangle\},
\end{split}
\end{equation}

\noindent It was shown in the paper {\rm \cite{Ir3}} (see Lemma 5 and Lemma 6) that the set of homology classes of $ H_n\left( K^H , K^H_{n-1}; \F\right)$ corresponded to the set of simplices 
$C^{\pi}_n(H)$ forms a basis of $ H_n\left( K^H , K^H_{n-1}; \F\right).$ Hence, the cardinality $|C^{\pi}_n(H)|$ doesn't depend on the choice of order $\pi$ on $\langle H \rangle.$

Let us fix any element $\langle w \rangle \in \langle H \rangle.$ Since $|C^{\pi}_n(H)|$ doesn't depend on the choice of $\pi,$ we can assume that $\langle w \rangle$ is the minimal element in $\pi$, i.e., $\pi(1) = \langle w \rangle.$ 
Let put 
\begin{equation}\label{eqC7}
\begin{split}
D^{\pi}_n(H; w) \MYdef \{\Delta = (\langle w \rangle,\langle w_{i_1}\rangle, \ldots ,\langle w_{i_n}\rangle) \in \langle H \rangle^{\times (n+1)}_{=0} \mid \\
1< i_1< \ldots < i_n\}.
\end{split}
\end{equation}
\noindent It is clear that 
\begin{equation}\label{eqC8}
D_n(H; w) \MYdef  |D^{\pi}_n(H; w)| 
\end{equation}
\noindent doesn't depend on the choice of $\pi.$

\begin{rem}\label{remC1}
From our definition, we have
\begin{equation}\label{remCC1}
D_n(H; w) = \dbinom{T-1}{n} - \dbinom{\langle H \rangle }{ \eta^{\bigstar}_n}^{\langle w \rangle}.
\end{equation}
\end{rem}

The map $\widehat{t} : C^{\pi}_n (H) \to \langle H \rangle^{\times n}$ defined on $ \Delta = (\langle w_{i_1}\rangle, \ldots ,\langle w_{i_{n+1}}\rangle) \in C^{\pi}_n (H) $ by the formula

$$\widehat{t} (\Delta) = (\langle w_{i_1}\rangle, \ldots , \langle w_{i_{t(\Delta)}-1}\rangle , \widehat{\langle w_{i_{t(\Delta)}}\rangle}, \langle w_{i_{t(\Delta)}+1}\rangle, \ldots ,
\langle w_{i_{n+1}}\rangle), $$

\noindent is a monomorphism. Note that for all $\Delta$ from the set
\begin{equation}\label{eqC9}
\begin{split}
C^{\pi}_{n; \ne 0}(H ; w) \MYdef \{\Delta = (\langle w_{i_1}\rangle, \ldots ,\langle w_{i_{n+1}}\rangle) \in C^{\pi}_n(H) \mid \\
span \ \langle w_{i_1}, \ldots , w_{i_{n+1}} \rangle  = n , \mbox{and} \, \langle w \rangle \notin span \ \langle w_{i_1}, \ldots , w_{i_{n+1}} \rangle \},
\end{split}
\end{equation}
\noindent we have
$$\widehat{t}(\Delta) \in \langle H \rangle^{\times n}_{\ne 0}.$$

\begin{lm}\label{lmC1}
The set of homology classes of $ H_n\left( K^H , K^H_{n-1}; \mathbb{Z}_2 \right)$ corresponded to the set of simplices $C^{\pi}_{n; \ne 0}(H; w)\cup D^{\pi}_n(H; w)$ generates $ H_n\left( K^H , K^H_{n-1}; \mathbb{Z}_2 \right).$
\end{lm}
\proof  Since the set of homology classes corresponded to the set of simplices $C^{\pi}_n(H)$ forms a basis of $ H_n\left( K^H , K^H_{n-1}; \mathbb{Z}_2 \right),$ it is sufficient to show that any simplex 
$\Delta = (\langle w_{i_1}\rangle, \ldots ,\langle w_{i_{n+1}}\rangle) \in C^{\pi}_n(H),$ such that 

$$\begin{array}{l} 
dim \ span \ \langle w_{i_1}, \ldots , w_{i_{n+1}} \rangle < n ,  \mbox{and} \   \langle w \rangle \notin \{ \langle w_{i_1}\rangle, \ldots , \langle w_{i_{n+1}}\rangle \},\\
 \mbox{or}\\
dim \ span \ \langle w_{i_1}, \ldots , w_{i_{n+1}} \rangle = n ,  \mbox{and} \  \langle w \rangle \in span \ \langle w_{i_1}, \ldots , w_{i_{n+1}} \rangle,\\
\mbox{but} \   \langle w \rangle \notin \{ \langle w_{i_1}\rangle, \ldots , \langle w_{i_{n+1}}\rangle \},\\
\end{array} $$
\noindent as a chain belongs to $span \langle C^{\pi}_{n; \ne 0}(H; w),  D^{\pi}_n(H; w), \delta_{n+1}(C_{n+1}(K^H; \mathbb{Z}_2))\rangle.$ In both cases, we have
$$\Delta_{n+1} \MYdef (\langle w \rangle, \langle w_{i_1}\rangle, \ldots ,\langle w_{i_{n+1}}\rangle) \in C_{n+1}(K^H; \mathbb{Z}_2),$$
\noindent and
$$\delta_{n+1}(\Delta_{n+1})= \Delta + \sum_{k=1}^{n+1} (\langle w \rangle, \langle w_{i_1}\rangle, \ldots ,\widehat{\langle w_{i_k}\rangle}, \ldots ,\langle w_{i_{n+1}}\rangle).$$

\noindent Hence, $\Delta \in span \langle  D^{\pi}_n(H; w), \delta_{n+1}(C_{n+1}(K^H; \mathbb{Z}_2))\rangle.$

\begin{flushright} {\sc Q.E.D.} \end{flushright}

\begin{cor} \label{corC1}
We have
\begin{equation}\label{eqC10}
|C^{\pi}_n(H)| \leq |C^{\pi}_{n; \ne 0}(H; w)| +  D_n(H; w).
\end{equation}
\end{cor}
\bigskip

We fix any $\langle w \rangle \in \langle H \rangle,$ and for any $\langle u \rangle \in \langle H \rangle,$ and order $\pi: [T] \to \langle H \rangle,$ such that $\pi (1) = \langle w \rangle$, we put

$$\langle H^{\pi}_{<u} \rangle \MYdef  \{\langle h \rangle \in \langle H \rangle \mid \langle h \rangle <_{\pi} \langle u \rangle \}.$$

\noindent For any subset $\langle U \rangle = \{\langle u_{i_{n-k}} \rangle, \ldots , \langle u_{i_1} \rangle \} \subset \langle H \rangle \setminus \langle w \rangle,$ 
$0 \leq k \leq n,$ such that  
\begin{equation}\label{eqC11}
\begin{split}
\langle u_{i_{n-k}} \rangle <_{\pi} \cdots  <_{\pi} \langle u_{i_1} \rangle , \\
dim \ span \langle U\rangle  \MYdef dim \ span \langle u_{i_{n-k}} , \ldots , u_{i_1} \rangle = n-k, \\
\mbox{and} \quad \langle w \rangle \notin \ span \langle U \rangle,
\end{split}
\end{equation}
\noindent we define the set $C^{\pi}_{n-k}(H^{\pi}_U; U) \subset \widehat{t} (C^{\pi}_{n; \ne 0}(H; w)),$ where $\langle H^{\pi}_U \rangle \MYdef \langle H^{\pi}_{<u_{i_{n-k}}}\rangle,$ in the following way
\begin{equation}\label{eqC12}
\begin{split}
C^{\pi}_{n-k}(H^{\pi}_U; U) \MYdef  \{\gamma \in \langle H\rangle^{\times n}_{\ne 0} \mid \\
\exists \Delta =(\langle w_{j_1} \rangle, \ldots , \langle w_{j_{k+1}} \rangle, \langle u_{i_{n-k}} \rangle, \ldots , \langle u_{i_1} \rangle) \in C^{\pi}_{n; \ne 0}(H; w) \, \mbox{s.t.} \, \widehat{t}(\Delta) = \gamma \} .
\end{split}
\end{equation}

\begin{rem}\label{remeqC12}
Note that from the conditions $ \Delta \in C^{\pi}_{n; \ne 0}(H; w)$ and $dim \ span \langle U\rangle = n-k,$ follows that $\langle w(\Delta) \rangle <_{\pi} \langle u_{i_{n-k}} \rangle <_{\pi} \cdots  <_{\pi} \langle u_{i_1} \rangle.$
\end{rem}
 
\begin{df}\label{dfCC1}
Let a subset $\langle U \rangle \subset \langle H\rangle \setminus \langle w \rangle$ satisfy the conditions (\ref{eqC11}) and $\langle U^{\pi}_s \rangle  \MYdef \{\langle u_{i_s} \rangle, \ldots , \langle u_{i_1} \rangle\},$ $1\leq s \leq n-k.$
We say that $\langle U \rangle$ is a $\pi$-closed $(n-k)$-subset of  $\langle H \rangle$ iff
\begin{equation}\label{eqC13}
span \langle U^{\pi}_s \rangle \cap \langle H^{\pi}_{U^{\pi}_s}\rangle = \O, \quad 1\leq s \leq n-k.
\end{equation}
\end{df}
\begin{df}\label{dfCC2}
We say that $\langle U \rangle$ is a $\pi$-boundary $(n-k)$-subset of  $\langle H \rangle$ iff
\begin{align}
& span \langle U^{\pi}_s \rangle \cap \langle H^{\pi}_{U^{\pi}_s}\rangle = \O, \quad 1\leq s \leq n-k-1, \label{eqCC14} \\
& span \langle U \rangle \cap \langle H^{\pi}_U\rangle \ne \O. \label{eqCC15}
\end{align}
\end{df}

Let $B(\langle H \rangle; n-k; \pi)$ denotes the set of all $\pi$-boundary $(n-k)$-subsets of  $\langle H \rangle,$ $0\leq k \leq n-2.$ Note that for $k=n-1,$ $B(\langle H \rangle; 1; \pi)= \O.$

\begin{lm}\label{lmC2}
The set $\widehat{t} (C^{\pi}_{n; \ne 0}(H; w)) \subset \langle H \rangle ^{\times n}_{\ne 0}$ can be expressed as disjoint union of $ C^{\pi}_{n-k}(H^{\pi}_U; U)$ as follows:
\begin{equation}\label{eqC14}
\widehat{t} (C^{\pi}_{n; \ne 0}(H; w)) = \bigsqcup_{k=0}^{n-2} \quad \bigsqcup_{\langle U \rangle \in B(\langle H \rangle; n-k; \pi) } C^{\pi}_{n-k}(H^{\pi}_U; U) .
\end{equation}
\end{lm}
\proof Let us assume that there exists $\gamma \in C^{\pi}_{n-k}(H^{\pi}_U; U) \cap C^{\pi}_{n-l}(H^{\pi}_V; V) \subset \langle H\rangle^{\times n}_{\ne 0}$ for some two different $\pi$-boundary subsets $\langle U \rangle$ and $\langle V \rangle.$ 
Then
\begin{align*}
&\widehat{t}^{-1}(\gamma)=\Delta = \\
&(\langle w_{j_1} \rangle, \ldots , \langle w_{j_{k+1}} \rangle, \langle u_{i_{n-k}} \rangle, \ldots , \langle u_{i_1} \rangle) = (\langle w_{r_1} \rangle, \ldots , \langle w_{r_{l+1}} \rangle, \langle v_{i_{n-l}} \rangle, \ldots , \langle v_{i_1} \rangle).
\end{align*}
Since $\langle U \rangle$ and $\langle V \rangle$ are different, then $k\ne l.$ Let $k<l.$ Then $\langle V \rangle \subset \langle U \rangle,$ $\langle v_{i_1} \rangle = \langle u_{i_1} \rangle, \ldots, \langle v_{i_{n-l}} \rangle = \langle u_{i_{n-l}} \rangle.$
Since $\langle U \rangle$ and $\langle V \rangle$ are $\pi$-boundary subsets, then
$$\O \ne span \langle V \rangle \cap \langle H^{\pi}_V \rangle = span \langle U^{\pi}_{n-l} \rangle \cap \langle H^{\pi}_{U^{\pi}_{n-l}} \rangle =\O.$$
This contradicts to our assumption. Hence, for any two different $\pi$-boundary subsets $\langle U \rangle$ and $\langle V \rangle,$ we have 
$$C^{\pi}_{n-k}(H^{\pi}_U; U) \cap C^{\pi}_{n-l}(H^{\pi}_V; V) =\O.$$

For any $\Delta \in C^{\pi}_{n; \ne 0}(H; w),$ $\Delta = (\langle w_{i_1} \rangle, \ldots, \langle w_{i_{n+1}} \rangle ),$ let
$$M(\Delta) \MYdef \{ t\in [n] \mid span \langle w_{i_{t+1}}, \ldots, w_{i_{n+1}} \rangle \cap \langle H^{\pi}_{<w_{i_{t+1}}} \rangle \ne \O \},$$
\noindent and
$$ m(\Delta) \MYdef \max_{t\in M(\Delta)} t. $$
\noindent Then the set $\langle U(\Delta) \rangle \MYdef \{\langle w_{i_{m(\Delta)+1}} \rangle, \ldots, \langle w_{i_{n+1}} \rangle\}$ is a $\pi$-boundary $(n-m(\Delta)+1)$-subset. Since $t(\Delta)\leq m(\Delta),$ then
$$\widehat{t} (\Delta) \in C^{\pi}_{n-m(\Delta)+1}(H^{\pi}_{U(\Delta)}; U(\Delta)).$$

\begin{flushright} {\sc Q.E.D.} \end{flushright}

Let denote by $\langle H ; w\notin\rangle^{\times s}_{\ne 0}$ the set
\begin{equation}\label{eqC15}
\langle H ; w\notin\rangle^{\times s}_{\ne 0} \MYdef \{U \in \langle H \rangle^{\times s}_{\ne 0} \mid w \notin span \langle U \rangle \}.
\end{equation}

\bigskip

\begin{thm}\label{thmC1}
Let $\langle H \rangle \subset {\bf RP}^n$ be a finite subset,  $|\langle H \rangle | = T.$  Then for any $\langle w \rangle \in  \langle H \rangle ,$ the following inequality is true  
\begin{equation}\label{eqC16}
\begin{split}
\rank \, H_n\left( K^H , K^H_{n-1}; \F\right) \leq \frac{1}{(T-1)!}\sum_{k=0}^{n-2} \quad \sum_{U\in \langle H ; w\notin\rangle^{\times (n-k)}_{\ne 0}} \frac{1}{q_{n-k-1}^U\cdots q_1^U} \times \\
\times \sum_{d=k+3}^{T-q_{n-k-1}^U} \frac{(T-d)!(T-q_{n-k}^U-1)!}{(T-q_{n-k-1}^U-d)!} \cdot A_{n-k}(U; d) {d-2 \choose k} + D_n(H; w),
\end{split}
\end{equation}
\noindent where
\begin{equation*}
\begin{split}
A_{n-k}(U; d) \MYdef \left[ {T-q_{n-k-1}^U -2 \choose q_{n-k}^U - q_{n-k-1}^U-1} - {T-q_{n-k-1}^U -d \choose q_{n-k}^U - q_{n-k-1}^U-1}\right] \times \\
\\
\times (q_{n-k}^U - q_{n-k-1}^U-1)! .
\end{split}
\end{equation*}
\end{thm}

\proof Let denote by $\Gamma$ the set of all orders on the set $\langle H \rangle$ such that $\tilde \gamma (1) = \langle w \rangle,$ $\forall \tilde \gamma \in \Gamma.$ For a fixed order $\pi$ (see (\ref{eq23})) such that $\pi (1) = \langle w \rangle,$ any 
$\tilde \gamma \in \Gamma$ can be expressed as composition
$$\tilde \gamma = \pi \circ \gamma : [T] \to [T] \to \langle H \rangle$$
\noindent $\pi$ with a permutation $\gamma \in Sym([T]; 1) \subset Sym([T])$ such that $\gamma(1) =1.$
Let $p$ be the uniform probabilty distrubution on $\Gamma \cong Sym([T];1)$ :
\begin{equation*}\label{eqC17}
p(\gamma) =\frac{1}{(T-1)!} , \quad \gamma \in Sym([T];1).
\end{equation*}
For any collection $U=(\langle u_{i_{n-k}} \rangle, \ldots , \langle u_{i_1}) \rangle \in H^{\times (n-k)}_{\ne 0}$ such that $\langle w \rangle \notin span \langle U \rangle,$ we define the random function $I_{n-k}^U : Sym([T]; 1) \to \R$ by the formula:
\begin{equation}\label{eqC18}
I_{n-k}^U(\gamma) \MYdef \left\{\begin{array}{ll}
\dbinom{(\pi \circ \gamma)^{-1}(\langle u_{i_{n-k}}\rangle) -2}{  k}, &\mbox{if} \,\,   U\in B(\langle H \rangle; n-k; \pi\circ \gamma); \\  
    & \\
    &\\
\qquad \qquad \quad 0 & \mbox{in all other cases.}
\end{array}\right. 
\end{equation}

From definition (\ref{eqC12}) follows that for any $(n-k)$-subset $\langle U \rangle \subset \langle H \rangle \setminus \langle w \rangle$ and order $\pi:[T] \to \langle H \rangle$ satisfying to (\ref{eqC11}), we have:
\begin{equation}\label{eqC19}
\left| C^{\pi}_{n-k}(H^{\pi}_U; U) \right| < \dbinom{\left|H^{\pi}_U\right|}{k} .
\end{equation}
Let
\begin{equation}\label{eqC20}
I_{n-k} \MYdef \sum_{U\in \langle H ; w\notin\rangle^{\times (n-k)}_{\ne 0}} I_{n-k}^U : Sym([T]; 1) \to \R .
\end{equation}
Then from Corollary \ref{corC1}, (\ref{eqC19}), and Lemma \ref{lmC2} for any $\gamma \in Sym([T]; 1),$ we have
\begin{equation}\label{eqC21}
\rank \, H_n\left( K^H , K^H_{n-1}; \F\right) = \left|C_n^{\pi \circ \gamma}(H)\right| \leq \sum_{k=0}^{n-2}I_{n-k}(\gamma) + D_n(H; w).
\end{equation}
Hence, the inequality (\ref{eqC21}) holds if we change the right hand side of  (\ref{eqC21}) by its expectation:
\begin{equation}\label{eqC22}
\rank \, H_n\left( K^H , K^H_{n-1}; \F\right)  \leq \sum_{k=0}^{n-2}\mathbb{E}[I_{n-k}] + D_n(H; w).
\end{equation}
Let $U=(\langle u_{i_{n-k}} \rangle, \ldots , \langle u_{i_1}) \rangle) \in H^{\times (n-k)}_{\ne 0}$ such that $\langle w \rangle \notin span \langle U \rangle$ and $q_{n-k}^U > q_{n-k-1}^U+1.$ 
We don't need to consider the case $q_{n-k}^U = q_{n-k-1}^U+1,$ because for any $\gamma \in Sym([T]; 1)$ such $U$ cannot be a $\pi \circ \gamma$-boundary subset of $\langle H\rangle.$

Let us calculate the number of permutations $\gamma \in Sym([T]; 1)$ such that
\begin{equation}\label{eqC23}
I_{n-k}^U(\gamma) = \dbinom{d-2}{k},
\end{equation}
for some $d$, $k+3 \leq d \leq T - q_{n-k-1}^U.$
From (\ref{eqC18}), we have
$$(\pi \circ \gamma)^{-1}(\langle u_{i_{n-k}}\rangle) =d.$$
\noindent From the condition (\ref{eqCC14}) follows that positions of elements from $L_{n-k-1}(U) \cap \langle H \rangle$ in the order $\pi \circ \gamma$ have to be chosen from the set $[d+1, T] = [d+1, d+2, \ldots, T],$
and if they are fixed, we have $A_{n-k}(U; d)$ possibilities for arrangement of elements from the set 
$$\{L_{n-k}(U)\cap\langle H \rangle\}\setminus \{\{\langle u_{n-k}\rangle \cup L_{n-k-1}(U)\} \cap \langle H \rangle\}$$
\noindent to fulfill the $\pi \circ \gamma$-boundary condition (\ref{eqCC15}).
The arrangement of the remaining elements from 
$$\langle H \rangle \setminus \{\{\langle w\rangle \cup L_{n-k}(U)\} \cap \langle H \rangle\}$$
\noindent will not affect the fact that $U$ is a $\pi \circ \gamma$-boundary subset of $\langle H \rangle.$ 
Since the collection $U_{n-k-1} \MYdef (\langle u_{i_{n-k-1}} \rangle, \ldots , \langle u_{i_1}) \rangle)$ has to be a $\pi \circ \gamma$-closed subset of $\langle H \rangle,$ then the element $\langle u_{i_{n-k-1}}\rangle$ has to be in the first place among any 
$q_{n-k-1}^U$ positions selected from the set $[d+1, T]$ for arrangement of the set $L_{n-k-1}(U) \cap \langle H \rangle,$ while $q_{n-k-2}^U$ elements from $L_{n-k-2}(U) \cap \langle H \rangle$ can be located in any of the remained $q_{n-k-1}^U -1$ 
places. The arrangement of the elements from
$$\{L_{n-k-1}(U)\cap\langle H \rangle\}\setminus \{\{\langle u_{i_{n-k-1}}\rangle \cup L_{n-k-2}(U)\} \cap \langle H \rangle\}$$
\noindent in $q_{n-k-1}^U - q_{n-k-2}^U -1$ places, left after choosing $q_{n-k-2}^U +1$ places for arrangement of the set $L_{n-k-2}(U) \cap \langle H \rangle$ and $\langle u_{i_{n-k-1}}\rangle,$ doesn't affect the fulfillment of $\pi \circ \gamma$-closeness condition for $U_{n-k-1}.$

Continuing the same way, we get that the number $N(U; d)$ of permutations $\gamma \in Sym([T]; 1),$ such that $(\pi \circ \gamma)^{-1}(\langle u_{i_{n-k}}\rangle) =d$ and $U$ is a $\pi \circ \gamma$- boundary subset of $\langle H \rangle,$ is

$$\begin{array}{l} 
N(U;d) = (T-q_{n-k}^U -1)! A_{n-k}(U; d) \dbinom{T-d}{q_{n-k-1}^U} \times \\
\\
\times \dbinom{q_{n-k-1}^U -1}{q_{n-k-2}^U} (q_{n-k-1}^U - q_{n-k-2}^U-1)! \times \cdots \times \dbinom{q_l^U -1}{q_{l-1}^U} (q_l^U - q_{l-1}^U -1)! \times \\
\\
\cdots \times \dbinom{q_2^U -1}{q_1^U} (q_2^U - q_1^U-1)! = \frac{(T-q_{n-k}^U -1)! A_{n-k}(U; d) (T-d)!}{(T-q_{n-k-1}^U -d)!} \frac{1}{q_{n-k-1}^Uq_{n-k-2}^U\cdots q_1^U} .  
\end{array}$$
\noindent Hence,
\begin{equation}\label{eqC24}
\begin{split}
\mathbb{E}[I_{n-k}^U] = \frac{1}{(T-1)!}  \frac{1}{q_{n-k-1}^U \cdots q_1^U} \times \\
\\
\times \sum_{d=k+3}^{T-q_{n-k-1}^U} \frac{(T-d)!(T-q_{n-k}^U -1)! }{(T-q_{n-k-1}^U -d)!} A_{n-k}(U; d) \dbinom{d-2}{k}
\end{split}
\end{equation}
The Theorem follows from (\ref{eqC22}), (\ref{eqC24}), and additivity of expectation.

\begin{flushright} {\sc Q.E.D.} \end{flushright}

Taking into account the inequality (\ref{eq22}) and definition   (\ref{eq26}), we return to elaborating of $\eta^{\bigstar}_n (\langle H \rangle)$ for $\langle H \rangle = \langle E_n \rangle \subset {\bf RP}^n$ (see (\ref{eq1})).

We define $\delta_{n, k}$, $k= 1, \ldots, n+1,$ as
\begin{equation}\label{eq96}
\delta_{n, k} \MYdef \frac{|\langle E_n \rangle^{\times k} \setminus \langle E_n \rangle^{\times k}_{\ne 0}|}{|\langle E_n \rangle^{\times k}|}, \quad k=1, \ldots, n+1.
\end{equation} 
We choose subspaces 
$$V_k \subset {\bf R}^{n+1}, \quad dim \ V_k = k, \quad k= 1, \ldots, n+1,$$
such that the orthogonal projectors 
\begin{equation*}\label{eq97}
\begin{split}
P_k: {\bf R}^{n+1} = V_k^{\perp} \oplus V_k \to   V_k^{\perp} \oplus V_k  = {\bf R}^{n+1}, \\
P_k(v) = v_2 , \quad \forall v=v_1+v_2 \in  V_k^{\perp} \oplus V_k , \quad k=1, \ldots, n+1,
\end{split}
\end{equation*}
satisfy the following conditions:

\noindent for any $k$ linear independent vectors $w_{i_1}, \ldots, w_{i_k} \in E_n,$ the vectors 
$$w^k_{i_s} \MYdef P_k(w_{i_s}), \quad s=1, \ldots, k,$$
are linear independent as well.

Let $E_{n, k}$ denote the set
\begin{equation}\label{eq98}
E_{n, k} \MYdef P_k(E_n).
\end{equation}

For $W^{k+1} = (\langle w_{i_1}^{k+1} \rangle, \ldots , \langle w_{i_k}^{k+1}\rangle) \in \langle E_{n, k+1}\rangle^{\times k}_{\ne 0},$ we use the following notations:
\begin{equation}\label{eq99}
\begin{split}
L(W^{k+1}) \MYdef  span \   \langle w_{i_1}^{k+1}, \ldots , w_{i_k}^{k+1} \rangle \subset V_{k+1} = {\bf R}^{k+1}; \\
q_k^{W^{k+1}} \MYdef |L(W^{k+1}) \cap E_{n, k+1}|; \\
E^m_{n, k+1} \MYdef  \left\{ W^{k+1} \in \langle E_{n, k+1} \rangle^{\times k}_{\ne 0} \mid q_k^{W^{k+1}} = k+m \right\}; \\
\gamma^m_{k+1} \MYdef \frac{|E^m_{n, k+1}|}{|\langle E_{n, k+1}\rangle^{\times k}_{\ne 0}|}; \\
B_{k+1} \MYdef \left\{ (\langle w_{i_1}^{k+1} \rangle, \ldots, \langle w_{i_{k+1}}^{k+1} \rangle ) \in \langle E_{n, k+1}\rangle^{\times (k+1)} \setminus \langle E_{n, k+1}\rangle^{\times (k+1)}_{\ne 0} | \right. \\
\left. (\langle w_{i_1}^{k+1} \rangle, \ldots, \langle w_{i_{k}}^{k+1} \rangle ) \in \langle E_{n, k+1}\rangle^{\times k}_{\ne 0} \right\}; \\
\epsilon_{k+1} \MYdef \frac{|B_{k+1}|}{|\langle E_{n, k+1}\rangle^{\times (k+1)}|}.
\end{split}
\end{equation}

Note that 
\begin{equation}\label{eq100}
\delta_{n, s} = \frac{|\langle E_{n, k+1} \rangle^{\times s} \setminus \langle E_{n, k+1} \rangle^{\times s}_{\ne 0}|}{|\langle E_{n, k+1} \rangle^{\times s}|}, \quad s=1, \ldots, k+1.
\end{equation} 

Since
\begin{multline*}
\langle E_{n, k+1}\rangle^{\times (k+1)} \setminus \langle E_{n, k+1}\rangle^{\times (k+1)}_{\ne 0} = \\
=\left\{\langle E_{n, k+1}\rangle^{\times (k+1)} \setminus \langle E_{n, k+1}\rangle^{\times (k+1)}_{\ne 0} \setminus B_{k+1} \right\} \cup B_{k+1},
\end{multline*}
then
\begin{equation}\label{eq101}
\delta_{n, k+1} = \delta_{n, k} + \epsilon_{k+1}.
\end{equation}
From definition (\ref{eq99}) we have
\begin{equation}\label{eq102}
|B_{k+1}| = \sum_{m=1}^{2^{k-1} -k} |E_{n, k+1}^m| m ,
\end{equation} 

\begin{equation}\label{eq103}
|E_{n, k+1}^m| =\gamma_{k+1}^m |\langle E_{n, k+1} \rangle^{\times k}_{\ne 0}|, 
\end{equation} 
and from (\ref{eq100})
\begin{equation}\label{eq104}
|\langle E_{n, k+1} \rangle^{\times k}_{\ne 0}| = (1- \delta_{n, k})  |\langle E_{n, k+1} \rangle^{\times k}|. 
\end{equation}
Hence, from (\ref{eq102}), (\ref{eq103}), (\ref{eq104}) we get
\begin{equation}\label{eq105}
\epsilon_{k+1} = (1- \delta_{n, k}) \sum_{m=1}^{2^{k-1} -k}  \gamma_{k+1}^m \frac{m}{2^n-k}. 
\end{equation}

\bigskip

\begin{df}\label{defC1}
We say that a vector $w\in \R^{n+1}$ is in general position to the sets $E_n =P_{n+1}(E_n), P_n(E_n), \ldots, P_2(E_n)$  iff for any $k, k=2, \ldots, n+1,$ and vectors $w_{i_1}, \ldots, w_{i_{k-1}} \in E_n \subset \R^{n+1},$ the vector $w^k \MYdef P_k(w)$ doesn't belong to the linear span of vectors $w_{i_1}^k, \ldots, w_{i_{k-1}}^k:$
$$w^k \notin span \langle w_{i_1}^k, \ldots, w_{i_{k-1}}^k \rangle .$$
\end{df}

Next, we are going to apply the Theorem \ref{thmC1} to the sets
$$\langle Z_k\rangle \MYdef \langle E_{n,k}\rangle \cup \langle w^k \rangle \subset {\bf RP}^{k-1}, \quad k=2, \ldots, n+1,$$
\noindent where $w\in \R^{n+1}$ is in general position to the sets $P_k(E_n),$ $k, k=2, \ldots, n+1.$
We express the right hand side of the inequality (\ref{eqC16}) as the sum of four summands
\begin{equation}\label{eqC25}
\rank \, H_{k-1}\left( K^{Z_k} , K^{Z_k}_{k-2}; \F\right)  \leq  S_0^k +S_{\leq}^k + S_>^k + D_{k-1}(Z_k; w^k),
\end{equation}
where
\begin{equation}\label{eqC26}
\begin{split}
S_0^k \MYdef \frac{1}{(2^n)!} \sum_{U\in \langle Z_k ; w^k\notin\rangle^{\times (k-1)}_{\ne 0}} \frac{1}{q_{k-2}^U\cdots q_1^U} \times \\
\\
\times \sum_{d=3}^{2^n +1 -q_{k-2}^U} \frac{(2^n+1-d)!(2^n-q_{k-1}^U)!}{(2^n +1-q_{k-2}^U-d)!} \cdot A_{k-1}(U; d);
\end{split}
\end{equation}
\\
\begin{equation}\label{eqC27}
\begin{split}
S_{\leq}^k \MYdef \frac{1}{(2^n)!} \sum_{t=1}^{k-3}\sum_{U\in \langle Z_k ; w^k\notin\rangle^{\times (k-t-1)}_{\ne 0}} \frac{1}{q_{k-t-2}^U\cdots q_1^U} \times \\
\\
\times \sum_{d=t+3}^{q_{k-t-2}^U+2} \frac{(2^n+1-d)!(2^n-q_{k-t-1}^U)!}{(2^n +1-q_{k-t-2}^U-d)!} \cdot A_{k-t-1}(U; d) \dbinom{d-2}{t};
\end{split}
\end{equation}
\\
\begin{equation}\label{eqC28}
\begin{split}
S_>^k \MYdef \frac{1}{(2^n)!} \sum_{t=1}^{k-3}\sum_{U\in \langle Z_k ; w^k\notin\rangle^{\times (k-t-1)}_{\ne 0}} \frac{1}{q_{k-t-2}^U\cdots q_1^U} \times \\
\\
\sum_{d=\max(q_{k-t-2}^U+3; t+3)}^{2^n +1 -q_{k-t-2}^U} \frac{(2^n+1-d)!(2^n-q_{k-t-1}^U)!}{(2^n +1-q_{k-t-2}^U-d)!} \cdot A_{k-t-1}(U; d) \dbinom{d-2}{t},
\end{split}
\end{equation}
\noindent where
\begin{equation*}
\begin{split}
A_{k-t-1}(U; d) \MYdef \left[ {2^n-q_{k-t-2}^U -1 \choose q_{k-t-1}^U - q_{k-t-2}^U-1} - {2^n-q_{k-t-2}^U -d+1 \choose q_{k-t-1}^U - q_{k-t-2}^U-1}\right] \times \\
\\
\times (q_{k-t-1}^U - q_{k-t-2}^U-1)! .
\end{split}
\end{equation*}

\begin{lm}\label{lmC3}
For sufficiently large $n,$ the following inequality holds:
\begin{equation}\label{eqC29}
S_>^k + S_{\leq}^k \leq \left( 1 + o\left(\frac{n^2}{2^n}\right)\right)\frac{k-3}{2^n}\dbinom{2^n}{k-1}.
\end{equation}
\end{lm}
\proof From (\ref{eqC28}) and the inequality $q_{k-t-2}^U \geq k-t-2$  we have 
\begin{equation}\label{eqC29}
\begin{split}
S_>^k \leq \frac{1}{2^n} \sum_{t=1}^{k-3}\sum_{U\in \langle Z_k ; w^k\notin\rangle^{\times (k-t-1)}_{\ne 0}} \frac{1}{q_{k-t-2}^U\cdots q_1^U} \times \\
\\
\times \sum_{d=\max(q_{k-t-2}^U+3; t+3)}^{2^n +1 -q_{k-t-2}^U} \frac{(2^n-d+1) \cdots (2^n-d - q_{k-t-2}^U +2)}{(2^n -1) \cdot(2^n -q_{k-t-2}^U)} \dbinom{d-2}{t} \leq \\
\\
\leq \frac{1}{2^n} \sum_{t=1}^{k-3}\sum_{U\in \langle Z_k ; w^k\notin\rangle^{\times (k-t-1)}_{\ne 0}} \frac{1}{q_{k-t-2}^U\cdots q_1^U} \frac{(2^n-1) \cdots (2^n-t)}{t!} \times \\
\\
\times \sum_{d=\max(q_{k-t-2}^U+3; t+3)}^{2^n +1 -q_{k-t-2}^U}(1-\alpha_d)^{q_{k-t-2}^U}\alpha_d^t \leq \\
\\
\leq \frac{1}{2^n} \sum_{t=1}^{k-3}\sum_{U\in \langle Z_k ; w^k\notin\rangle^{\times (k-t-1)}_{\ne 0}} \frac{1}{q_{k-t-2}^U\cdots q_1^U} \frac{(2^n-1) \cdots (2^n-t)}{t!} \times \\
\\
\times \sum_{d=\max(q_{k-t-2}^U+3; t+3)}^{2^n +1 -q_{k-t-2}^U}(1-\alpha_d)^{k-t-2}\alpha_d^t .
\end{split}
\end{equation}
Here $\alpha_d = \frac{d-2}{2^n-1}.$ 

 From (\ref{eqC27}) we have
\begin{equation}\label{eqC30}
\begin{split}
S_{\leq}^k \leq \frac{1}{2^n} \sum_{t=1}^{k-3}\sum_{U\in \langle Z_k ; w^k\notin\rangle^{\times (k-t-1)}_{\ne 0}} \frac{1}{q_{k-t-2}^U\cdots q_1^U} \times \\
\\
\times \sum_{d= t+3}^{q_{k-t-2}^U+2} \frac{(2^n- q_{k-t-2}^U-1) \cdots (2^n - q_{k-t-2}^U -(d-2))}{(2^n -1) \cdot(2^n -(d-2))} \dbinom{d-2}{t} \leq \\
\\
\leq \frac{1}{2^n} \sum_{t=1}^{k-3}\sum_{U\in \langle Z_k ; w^k\notin\rangle^{\times (k-t-1)}_{\ne 0}} \frac{1}{q_{k-t-2}^U\cdots q_1^U} \frac{(2^n-1) \cdots (2^n-t)}{t!} \times \\
\\
\times \sum_{d= t+3}^{q_{k-t-2}^U+2} \left(1-\frac{q_{k-t-2}^U}{2^n-1}\right)^{d-2}\alpha_d^t \leq \\
\\
\leq \frac{1}{2^n} \sum_{t=1}^{k-3}\sum_{U\in \langle Z_k ; w^k\notin\rangle^{\times (k-t-1)}_{\ne 0}} \frac{1}{q_{k-t-2}^U\cdots q_1^U} \frac{(2^n-1) \cdots (2^n-t)}{t!} \times \\
\\
\times \sum_{d= t+3}^{q_{k-t-2}^U+2} (1- \alpha_d)^{k-t-2}\alpha_d^t.
\end{split}
\end{equation}
The last inequality in (\ref{eqC30}) follows from the inequality
$$\left(1- \frac{q}{2^n-1}\right)^{(2^n-1)x} \leq (1-x)^{k-t-2}$$
that holds for $x\in \left[0, \frac{q}{2^n-1}\right]$ and $q\geq k-t-2.$ 

From (\ref{eqC29}) and (\ref{eqC30}) we get
\begin{equation*}
\begin{split}
S_>^k + S_{\leq}^k \leq \frac{1}{2^n} \sum_{t=1}^{k-3} \frac{(2^n-1) \cdots (2^n-t)}{t!} \sum_{d=2}^{2^n+1} (1-\alpha_d)^{k-t-2}\alpha_d^t \times \\
\\
\times \sum_{U\in \langle Z_k ; w^k\notin\rangle^{\times (k-t-1)}_{\ne 0}} \frac{1}{q_{k-t-2}^U\cdots q_1^U}  \leq \\
\\
\leq \frac{1}{2^n} \sum_{t=1}^{k-3} \frac{(2^n-1) \cdots (2^n-t)}{t!} \left( \frac{\Gamma(k-t-1)\Gamma(t+1)}{\Gamma(k)} + \frac{n^2}{12(2^n-1)^2} \right) \times
\\
\times \frac{1}{(k-t-2)!} \sum_{U\in \langle Z_k ; w^k\notin\rangle^{\times (k-t-1)}_{\ne 0}} 1 \leq \left(1 + o\left( \frac{n^2}{2^n}\right)\right) \frac{k-3}{2^n} \dbinom{2^n}{k-1}.
\end{split}
\end{equation*}
Here we used the formula 853.21 from \cite{Dwi} that holds for $t\geq 1$ :
$$\int\limits_0^1 (1-x)^{k-t-2} x^t \, dx = \frac{\Gamma(k-t-1)\Gamma(t+1)}{\Gamma(k)} = \frac{(k-t-2)!t!}{(k-1)!}$$
\noindent and the midpoint rule for sum estimation by integral.

\begin{flushright} {\sc Q.E.D.} \end{flushright}

\begin{lm}\label{lmC4}
For any $n\geq 64,$ the following inequality holds:
\begin{equation}\label{eqC31}
S_0^k  \leq \left(  \frac{1}{2^n}\left( 1 + o\left(\frac{n^3}{2^n}\right) \right)+ \left(\delta_{n, k} -\delta_{n, k-1}\right)\left(1 - \frac{k-1}{2^n}\right)\right)\dbinom{2^n}{k-1}.
\end{equation}
\end{lm}
\proof   From (\ref{eqC26}) we have
\begin{equation}\label{eqC32}
\begin{split}
S_0^k \leq \frac{1}{2^n} \sum_{U\in \langle Z_k ; w^k\notin\rangle^{\times (k-1)}_{\ne 0}} \frac{1}{q_{k-2}^U\cdots q_1^U} \left(J(U, k; >) + J(U, k ; \leq)\right),
\end{split}
\end{equation}
where
\begin{equation}\label{eqC33}
\begin{split}
J(U, k; >) \MYdef \sum_{d=q_{k-2}^U+3}^{2^n +1 -q_{k-2}^U} \frac{(2^n-d+1) \cdots (2^n-d - q_{k-2}^U +2)}{(2^n -1) \cdots (2^n -q_{k-2}^U)},
\end{split}
\end{equation}
and
\begin{equation}\label{eqC34}
\begin{split}
J(U, k ; \leq) \MYdef \sum_{d=3}^{q_{k-2}^U+2}\Bigg( \frac{(2^n -q_{k-2}^U -1) \cdots (2^n-q_{k-2}^U -(d-2))}{(2^n-1) \cdots (2^n-(d-2))} - \\
\\
- \frac{(2^n -q_{k-1}^U ) \cdots (2^n-q_{k-1}^U -(d-3))}{(2^n-1) \cdots (2^n-(d-2))}\Bigg).
\end{split}
\end{equation}
We have
\begin{equation}\label{eqC35}
\begin{split}
\frac{1}{q_{k-2}^U\cdots q_1^U} J(U, k; >) \leq \frac{1}{q_{k-2}^U\cdots q_1^U} \sum_{d=q_{k-2}^U+3}^{2^n +1 -q_{k-2}^U} (1- \alpha_d)^{q_{k-2}^U} \leq \\
\\
\leq \frac{1}{(k-2)!} \frac{k-2}{q_{k-2}^U} \left( \frac{1}{q_{k-2}^U +1} + \frac{\left(q_{k-2}^U\right)^2}{24(2^n-1)^2}\right) \leq \\
\\
\leq \frac{1}{(k-1)!} \left( 1 + \frac{n^2 q_{k-2}^U}{(2^n-1)^2}\right) = \frac{1}{(k-1)!}\left( 1 + o\left(\frac{n^3}{2^n}\right)\right).
\end{split}
\end{equation}
Here $\alpha_d = \frac{d-2}{2^n-1}.$

\bigskip

We assert that for any $U\in \langle Z_k ; w^k\notin\rangle^{\times (k-1)}_{\ne 0}$ for $n\geq 64,$ the following inequality holds:
\begin{equation}\label{eqD37}
\begin{split}
\frac{1}{q_{k-2}^U\cdots q_1^U} J(U, k; \leq) \leq \frac{q_{k-1}^U - (k-1)}{(k-1)!}.
\end{split}
\end{equation}

For $U\in \langle Z_k ; w^k\notin\rangle^{\times (k-1)}_{\ne 0}$ such that $q_{k-1}^U \geq n^3,$ we have 
\begin{equation}\label{eqC37}
\begin{split}
\frac{1}{q_{k-1}^U - (k-1)}\frac{q_{k-1}^U - (k-1)}{q_{k-2}^U\cdots q_1^U} J(U, k; \leq) \leq \\
\\
\leq \frac{1}{q_{k-1}^U - (k-1)}\frac{1}{q_{k-2}^U}\frac{q_{k-1}^U - (k-1)}{(k-3)!}\cdot q_{k-2}^U \leq \\
\\
(q_{k-1}^U \geq n^3)\\
\\
\leq \frac{q_{k-1}^U - (k-1)}{(k-1)!}.
\end{split}
\end{equation}

For $U\in \langle Z_k ; w^k\notin\rangle^{\times (k-1)}_{\ne 0}$ such that $q_{k-1}^U < n^3,$ it can be proven by mathematical induction that for any $d,  \ 3\leq d \leq q_{k-2}^U +2,$ and $n \geq 64$ (for $n \geq 64$, we have $n^{10}<2^n-1$), the following inequality holds:
\begin{equation}\label{eqC38}
\begin{split}
\frac{(2^n -q_{k-2}^U -1) \cdots (2^n-q_{k-2}^U -(d-2))}{(2^n-1) \cdots (2^n-(d-2))} - \\ 
\\
- \frac{(2^n -q_{k-1}^U ) \cdots (2^n-q_{k-1}^U -(d-3))}{(2^n-1) \cdots (2^n-(d-2))} \leq \frac{d-2}{n\left(q_{k-2}^U\right)^2} .
\end{split}
\end{equation}
Hence, for $U\in \langle Z_k ; w^k\notin\rangle^{\times (k-1)}_{\ne 0}$ such that $q_{k-1}^U < n^3,$ we have:
\begin{equation}\label{eqC39}
\begin{split}
\frac{1}{q_{k-1}^U - (k-1)}\frac{q_{k-1}^U - (k-1)}{q_{k-2}^U\cdots q_1^U} J(U, k; \leq) \leq \\
\\
\leq 1 \cdot \frac{q_{k-1}^U - (k-1)}{(k-2)!} \frac{1}{n\left(q_{k-2}^U\right)^2} \sum_{d=3}^{q_{k-2}^U+2}(d-2) \leq \frac{q_{k-1}^U - (k-1)}{(k-1)!}.
\end{split}
\end{equation}
The assertion (\ref{eqD37}) follows from (\ref{eqC37}) and (\ref{eqC39}).
\\

It follows from (\ref{eqC32}), (\ref{eqC35}), and (\ref{eqD37}) that for any $n\geq 64,$ the following inequality holds:
\begin{equation}\label{eqC40}
\begin{split}
S_0^k \leq \frac{1}{2^n} \frac{1 + o\left(\frac{n^3}{2^n}\right)}{(k-1)!} \sum_{U\in \langle Z_k ; w^k\notin\rangle^{\times (k-1)}_{\ne 0}} 1  \quad + \\
\\
+ \frac{1}{2^n} \frac{1}{(k-1)!}\sum_{U\in \langle Z_k ; w^k\notin\rangle^{\times (k-1)}_{\ne 0}}\left(q_{k-1}^U -(k-1)\right).
\end{split}
\end{equation}
Taking into account (\ref{eq101}), we note that
\begin{equation}\label{eqC41}
\begin{split}
\frac{1}{2^n}\sum_{U\in \langle Z_k ; w^k\notin\rangle^{\times (k-1)}_{\ne 0}}\left(q_{k-1}^U -(k-1)\right) = \frac{1}{2^n} |B_k| =\\
\\
= \frac{\epsilon_k}{2^n} 2^n(2^n-1)\cdots(2^n-k+1) = (\delta_{n, k} -  \delta_{n, k-1})\dbinom{2^n-1}{k-1}(k-1)!.
\end{split}
\end{equation}

The Lemma follows from (\ref{eqC40}) and (\ref{eqC41}).
\begin{flushright} {\sc Q.E.D.} \end{flushright}

\begin{thm}\label{thmC2}
For any $n\geq 64,$ the following inequality holds for $k=1, \ldots, n:$
\begin{equation}\label{eqC42}
\begin{split}
\eta^{\bigstar}_k(\langle Z_{k+1} \rangle ) \geq \left[ 1- \delta_{n, k} - \frac{k-1}{2^n}\left( 1 + o\left(\frac{n^3}{2^n}\right) \right)-  \left(1 - \frac{k}{2^n}\right)\left(\delta_{n, k+1} -\delta_{n, k}\right)\right]\dbinom{2^n}{k}=\\
\\
= \dbinom{\langle E_n \rangle }{ \eta^{\bigstar}_k}^{\langle w \rangle} - \left[\frac{k-1}{2^n}\left( 1 + o\left(\frac{n^3}{2^n}\right) \right) + \left(1 - \frac{k}{2^n}\right)\left(\delta_{n, k+1} -\delta_{n, k}\right)\right]\dbinom{2^n}{k} .
\end{split}
\end{equation}
\end{thm}
\proof From the exact sequence (\ref{eqC1}), we get:
\begin{equation}\label{eqC43}
\eta^{\bigstar}_k(\langle Z_{k+1} \rangle ) = \dbinom{2^n}{k} - rank \, H_k\left( K^{Z_{k+1}} , K^{Z_{k+1}}_{k-1}; \F\right).
\end{equation}
From inequality (\ref{eqC25}), Lemma \ref{lmC3}, and Lemma \ref{lmC4}, for any $n\geq 64,$  we get for $k=1, \ldots, n:$
\begin{equation}\label{eqC44}
\begin{split}
 rank \, H_k\left( K^{Z_{k+1}} , K^{Z_{k+1}}_{k-1}; \F\right) \leq \\
\\
\leq \left[ \delta_{n, k} + \frac{k-1}{2^n}\left( 1 + o\left(\frac{n^3}{2^n}\right) \right)+  \left(1 - \frac{k}{2^n}\right)\left(\delta_{n, k+1} -\delta_{n, k}\right)\right]\dbinom{2^n}{k}
\end{split}
\end{equation}
The Theorem follows from (\ref{eqC43}), (\ref{eqC44}), and Remark \ref{remC1}.

\begin{flushright} {\sc Q.E.D.} \end{flushright}

\vskip 1cm

\section{Asymptotics of the number of singular $\{\pm 1\}$-matrices.}

\begin{lm}\label{lm3}
For $n\geq 64$  and $k=1, \ldots, n,$ we have
\begin{equation}\label{eq106}
\delta_{n, k+1} -\delta_{n, k} \leq \frac{k-1}{2^n}\left( 1 + o\left(\frac{n^3}{2^n}\right) \right) .
\end{equation}
\end{lm}
\proof  We define the probability distribution $p: \langle Z_{k+1} \rangle =\langle E_{n, k+1} \rangle \cup \left\{\langle w^{k+1} \rangle \right\} \to [0, 1]$ by the rule:
\begin{equation*}\label{eq108}
p\left(\langle w^{k+1} \rangle \right) = 1, \quad p\left(\langle w^{k+1}_i \rangle \right) = 0, \quad i=1, \ldots, 2^n.
\end{equation*}
Then from Theorem \ref{thm4}  we have
\begin{equation}\label{eq109}
\eta^{\bigstar}_k \left( \langle Z_{k+1} \rangle \right) = \sum_{W^{k+1} \in \langle E_{n, k+1}\rangle^{\times k}_{\ne 0}} \frac{1}{W^{k+1}[E_{n, k+1}]}.
\end{equation}

For any permutation $\sigma \in Sym[k],$ we define the map $\sigma^*: \langle E_{n, k+1}\rangle^{\times k}_{\ne 0} \to \langle E_{n, k+1}\rangle^{\times k}_{\ne 0}$ by the formula
\begin{multline*}\label{eq110}
 \sigma^*\left(W^{k+1}\right) = W^{k+1}_{\sigma} \MYdef (\langle w_{i_{\sigma(1)}}^{k+1} \rangle, \ldots , \langle w_{i_{\sigma(k)}}^{k+1}\rangle) \in \langle E_{n, k+1}\rangle^{\times k}_{\ne 0},\\
\forall \  W^{k+1} = (\langle w_{i_1}^{k+1} \rangle, \ldots , \langle w_{i_k}^{k+1}\rangle) \in \langle E_{n, k+1}\rangle^{\times k}_{\ne 0}.
\end{multline*}
Note that for the symmetrization of a combinatorial flag, defined by the formula
$$ Sym\left( W^{k+1}\right) \MYdef \sum_{\sigma \in Sym[k]} \frac{1}{W^{k+1}_{\sigma}[E_{n, k+1}]},$$
for $W^{k+1}\in E_{n, k+1}^m,$ we have
\begin{equation}\label{eq111}
Sym\left( W^{k+1}\right)  \leq \frac{k}{k+m}.
\end{equation}
Hence,
\begin{equation}\label{eq112}
\begin{split}
\sum_{W^{k+1} \in E_{n, k+1}^m} \frac{1}{W^{k+1}[E_{n, k+1}]} \leq \gamma^m_{k+1}(1-\delta_{n, k})\frac{k}{k+m} \binom {2^n}{k} \quad \mbox{and} \\
\eta^{\bigstar}_k \left( \langle Z_{k+1} \rangle \right)\leq (1-\delta_{n, k}) \binom {2^n}{k}\sum_{m=0}^{2^{k-1} -k}\gamma^m_{k+1}\frac{k}{k+m}.
\end{split}
\end{equation}
Combaining (\ref{eq112}) with inequality (\ref{eqC42}), for $n\geq 64,$ we get 
\begin{equation}\label{eq116}
\begin{split}
(1-\delta_{n, k}) \binom {2^n}{k}\sum_{m=0}^{2^{k-1} -k}\gamma^m_{k+1}\frac{k}{k+m} \geq \eta^{\bigstar}_k(\langle Z_{k+1} \rangle ) \geq \\
\\
\geq \left[ 1- \delta_{n, k} - \frac{k-1}{2^n}\left( 1 + o\left(\frac{n^3}{2^n}\right) \right)-  \left(1 - \frac{k}{2^n}\right)\left(\delta_{n, k+1} -\delta_{n, k}\right)\right]\dbinom{2^n}{k},
\end{split}
\end{equation}
or
\begin{equation}\label{eqC116}
\begin{split}
\frac{1}{1- \delta_{n, k}} \left[\left( 1 + o\left(\frac{n^3}{2^n}\right) \right) \frac{k-1}{2^n} + \left(1 - \frac{k}{2^n}\right)\left(\delta_{n, k+1} -\delta_{n, k}\right)\right] \geq \\
\\
\geq 1- \gamma_{k+1}^0  - \sum_{m=1}^{2^{k-1} -k}\gamma^m_{k+1}\frac{k}{k+m}.
\end{split}
\end{equation}
Taking into account the identity
\begin{equation*}\label{eq117}
1- \gamma^0_{k+1} = \sum_{m=1}^{2^{k-1}-k} \gamma^m_{k+1},
\end{equation*}
the inequality (\ref{eqC116}) may be expressed as follows
\begin{equation}\label{eq118}
\begin{split}
\frac{1}{1- \delta_{n, k}} \left[\left( 1 + o\left(\frac{n^3}{2^n}\right) \right) \frac{k-1}{2^n} +  \left(1 - \frac{k}{2^n}\right)\left(\delta_{n, k+1} -\delta_{n, k}\right)\right] \geq \\
\\
\geq \sum_{m=1}^{2^{k-1} -k}\gamma^m_{k+1}\frac{m}{k+m}.
\end{split}
\end{equation}

For $1 \leq k \leq n-1,$ the inequality
\begin{equation}\label{eqD121}
\frac{m}{k+m} \geq \frac{2m}{2^n-k}
\end{equation}
is fullfiled for all $m,$ $1\leq m \leq 2^{k-1} -k.$ 

For $k=n,$ the inequality (\ref{eqD121}) is true for all $m\leq 2^{n-1} - \frac{3}{2}n.$  It follows from Littlewood-Offord lemma in the form proven by P. Erd\" os {\rm \cite{Erd}} that for $n\geq 4,$
\begin{equation}\label{eqD122}
\gamma^m_{n+1} =0, \quad \forall \  m, \ \mbox{such that} \quad 2^{n-1} - \frac{3}{2}n < m < 2^{n-1}-n.
\end{equation}

From (\ref{eq118}), (\ref{eqD121}),  (\ref{eqD122}), (\ref{eq101}), and (\ref{eq105}), we have
\begin{equation}\label{eqC122}
\begin{split}
\frac{1}{1- \delta_{n, k}} \left[\left( 1 + o\left(\frac{n^3}{2^n}\right) \right) \frac{k-1}{2^n} + \left(1 - \frac{k}{2^n}\right)\left(\delta_{n, k+1} -\delta_{n, k}\right)\right] \geq \\
\\
\geq \sum_{m=1}^{2^{k-1} -k}\gamma^m_{k+1}\frac{2m}{2^n-k} = \frac{2(\delta_{n, k+1} -\delta_{n, k})}{1- \delta_{n, k}}
\end{split}
\end{equation}
Hence,
\begin{equation}\label{eqC123}
\left( 1 + o\left(\frac{n^3}{2^n}\right) \right) \frac{k-1}{2^n} \geq \left(1 + \frac{k}{2^n}\right) \left(\delta_{n, k+1} -\delta_{n, k}\right). 
\end{equation}
The inequality (\ref{eq106}) of Lemma \ref{lm3} follows from (\ref{eqC123}).

\begin{flushright} {\sc Q.E.D.} \end{flushright}

For ease of use of established terminology, we formulate an estimate for the cardinality of the set of singular $\{\pm 1\}$-matrices in terms of the probability ${\mathbb P}_n$  of singularity of random Bernoulli matrices.
Also, we can identify $\langle E_n \rangle$ with $E_n$.
\begin{thm}\label{thm5}
For $n\to\infty,$ we have
\begin{equation}\label{eq124}
{\mathbb P}_n \thicksim \frac{(n-1)^2}{2^{n-1}}.
\end{equation}
\end{thm}
\proof  By definition, we have
\begin{equation}\label{eq125}
\begin{split}
{\mathbb P}_{n+1} = \frac{\left|\left[\{E_n\}^{n+1}\setminus [E_n]^{\times (n+1)}\right] \cup \left[[E_n]^{\times (n+1)} \setminus [E_n]^{\times (n+1)}_{\ne 0} \right]\right|}{2^{n(n+1)}},  \ \mbox{i.e.,} \\
\\
{\mathbb P}_{n+1} = \frac{\left|\{E_n\}^{n+1}\setminus [E_n]^{\times (n+1)}\right|}{2^{n(n+1)}} + \delta_{n, n+1}\frac{\left| [E_n]^{\times (n+1)}\right|}{2^{n(n+1)}},
\end{split}
\end{equation}
\noindent where $\{E_n\}^{n+1}= \underbrace{E_n\times \cdots \times E_n}_{n+1}.$

Cardinality of the subset of matrices containing exactly two equal rows asymptotically plays the main role for estimation of $\left|\{E_n\}^{n+1}\setminus [E_n]^{\times (n+1)}\right|,$ i.e.,
\begin{equation}\label{eq126}
\frac{\left|\{E_n\}^{n+1}\setminus [E_n]^{\times (n+1)}\right|}{2^{n(n+1)}} = \frac{n(n+1)}{2^{n+1}}(1+o_n(1)).
\end{equation}

From Lemma \ref{lm3} we have
\begin{equation}\label{eq127}
\begin{split}
\delta_{n, n+1} \leq \delta_{n, n} + \left( 1 + o\left(\frac{n^3}{2^n}\right) \right) \frac{n-1}{2^n} \leq  
\\
\leq \delta_{n, n-1} + \left( 1 + o\left(\frac{n^3}{2^n}\right) \right) \left(\frac{n-2}{2^n} + \frac{n-1}{2^n}\right)  \leq \ldots\\ 
\\
\leq  \left( 1 + o\left(\frac{n^3}{2^n}\right) \right) \sum_{k=1}^{n-1}\frac{k}{2^n} = \left( 1 + o\left(\frac{n^3}{2^n}\right) \right) \frac{n(n-1)}{2^{n+1}}.
\end{split}
\end{equation}

We need to show that
\begin{equation}\label{eq129}
\delta_{n, n+1} \geq \frac{n(n-1)}{2^{n+1}}(1+o_n(1)).
\end{equation}

Let $R_{n-1}^{n+1} \subset  [E_{n-1}]^{\times (n+1)}$ be the subset of ordered collections $W = ( w_{i_1}, \ldots,  w_{i_{n+1}}) \in [E_{n-1}]^{\times (n+1)}$ such that the columns ${\bar 1}, Y_2, \ldots, Y_n$ of the matrix $M(W)$ with rows $ w_{i_1}, \ldots,  w_{i_{n+1}}$
$$M(W) =\begin{pmatrix}
w_{i_1}\\
\vdots \\
w_{i_{n+1}}
\end{pmatrix} =
({\bar 1}, Y_2, \ldots, Y_n)
$$
\noindent are not collinear, i.e., $Y_i \ne \pm Y_j,$ $\forall \ i\ne j.$

We can construct an ordered collection $W^{\prime} \in [E_n]^{\times (n+1)}$ by placing a column $\pm Y_i$, $i=2, \ldots, n$, in one of $n$ positions:
$$M(W^{\prime}) = ({\bar 1}, Y_2, \ldots, Y_i, \ldots, \pm Y_i, \ldots, Y_n).$$
\noindent Then the total number of $W^{\prime} \in [E_n]^{\times (n+1)}$ such that the matrix $M(W^{\prime})$ has exactly two equal up to sign columns is not less than
\begin{equation}\label{eq130}
\frac{2(n-1)n}{2}\left|R_{n-1}^{n+1}\right|= n(n-1)\left|R_{n-1}^{n+1}\right|.
\end{equation}
Since
\begin{equation}\label{eq131}
\begin{split}
\frac{\left|[E_{n-1}]^{\times (n+1)}\right| }{\left|[E_{n}]^{\times (n+1)}\right| } = \frac{1}{2^{n+1}}(1+o_n(1)) \quad \mbox{and} \\
\\
\left|[E_{n-1}]^{\times (n+1)}\right| = (1+o_n(1))\left|R_{n-1}^{n+1}\right|,
\end{split}
\end{equation}
then (\ref{eq129}) follows from (\ref{eq130}) and (\ref{eq131}). The Theorem follows from (\ref{eq125}), (\ref{eq126}), (\ref{eq127}), and (\ref{eq129}).

\begin{flushright} {\sc Q.E.D.} \end{flushright}

\vskip 1cm

\section{Asymptotics of the number of threshold functions.}

In this section we use notations from the previous section.
  
\begin{thm}\label{thm7}
Asymptotics of the number of threshold functions is equal to $2{2^n-1 \choose n}:$
\begin{equation}\label{eq15}
P(2, n) \thicksim 2 {2^n-1 \choose n}, \quad n\to \infty.
\end{equation}
\end{thm}

\proof  We write the inequality (\ref{eqC42}) of the Theorem \ref{thmC2} for $k=n$ taking into account the inequality (\ref{eq106}) of Lemma \ref{lm3}:
\begin{equation}\label{eqC132}
\begin{split}
\eta^{\bigstar}_n(\langle E_{n} \rangle \cup \left\{\langle w \rangle \right\}) \geq 
\\
\geq \left[ 1- \delta_{n, n} - \frac{n-1}{2^{n-1}}\left( 1 + o\left(\frac{n^3}{2^n}\right) \right) \right]\left[\dbinom{2^n-1}{n}+\dbinom{2^n-1}{n-1}\right] \geq\\
\\
\geq \left[ 1- \frac{(n-1)^2}{2^{n}}\left( 1 + o\left(\frac{n^3}{2^n}\right) \right) \right]\dbinom{2^n-1}{n} +
\\
+ \left[ 1- \delta_{n, n} - \frac{n-1}{2^{n-1}}\left( 1 + o\left(\frac{n^3}{2^n}\right) \right) \right]\dbinom{2^n-1}{n-1}.\\
\\
\end{split}
\end{equation}

From Theorem \ref{thm1} we have:
\begin{equation}\label{eq134}
\eta^{\bigstar}_n \left( \langle E_{n} \rangle\right) = \eta^{\bigstar}_n \left( \langle E_{n} \rangle \cup \left\{\langle w\rangle \right\}\right) - \eta^{\bigstar}_{n-1} \left( \langle E_{n} \rangle^{\perp w}\right).
\end{equation}
 From Theorem \ref{thm2} we have:
\begin{equation}\label{eq135}
\eta^{\bigstar}_{n-1} \left( \langle E_{n} \rangle^{\perp w}\right) \leq\dbinom{\langle E_n \setminus \{u\}\rangle^{\perp w} }{ \eta^{\bigstar}_{n-1}}^{\langle u^{\perp w} \rangle},
\end{equation}
\noindent where $u = {\bar 1} \in E_n.$

A summand $ \eta^{\bigstar}_{n-1} \left(\{\langle u^{\perp w} \rangle, \langle w_{i_1}^{\perp w} \rangle, \ldots, \langle w_{i_{n-1}}^{\perp w} \rangle\}\right)$ from the right side of (\ref{eq135}) is equal to $1$ iff
\begin{equation}\label{eq136}
\begin{split}
span \ \langle u^{\perp w}, w_{i_1}^{\perp w}, \ldots, w_{i_{n-1}}^{\perp w}\rangle = \langle w \rangle^{\perp} = {\bf R}^{n}, \  \mbox{or} \\
dim \ span \ \langle u,  w_{i_1}, \ldots, w_{i_{n-1}}\rangle = n.
\end{split}
\end{equation}
It follows from symmetry of $E_n,$  (\ref{eq136}), and definition of $\delta_{n, n}$ that the right side of (\ref{eq135}) is equal to $(1-\delta_{n, n})\binom{2^n-1}{n-1}$ and
\begin{equation}\label{eq137}
\eta^{\bigstar}_{n-1} \left( \langle E_{n} \rangle^{\perp w}\right) \leq (1-\delta_{n, n})\binom{2^n-1}{n-1}.
\end{equation}
Taking into account (\ref{eq22}), (\ref{eq26}), (\ref{eqC132}),  (\ref{eq134}), and (\ref{eq137}), we get a lower bound for $P(2, n):$
\begin{equation}\label{eq138}
P(2, n) \geq  2\left[ 1- \frac{n^2}{2^n}\left( 1 + o\left(\frac{n^3}{2^n}\right) \right) \right]\dbinom{2^n-1}{n}.
\end{equation}
The Theorem follows from the upper bound (\ref{eq2}) and the lower bound (\ref{eq138}).

\begin{flushright} {\sc Q.E.D.} \end{flushright}

\bigskip

\noindent Faculty of Mechanics and Mathematics, Lomonosov Moscow State University, Moscow, Russian Federation

\bigskip

{\it E-mail address}: irmatov@intsys.msu.ru

\end{document}